\documentclass[english,ruled]{article}
\usepackage[T1]{fontenc}
\usepackage[latin9]{inputenc}
\usepackage{verbatim}
\usepackage{algorithm2e}
\usepackage{amsmath}
\usepackage{amsthm}
\usepackage{amssymb}
\usepackage{graphicx}
\usepackage{xcolor}
\usepackage{multicol}
\usepackage{multirow}
\usepackage{geometry}
\usepackage{booktabs}
\usepackage{enumitem}
\usepackage{setspace}
\makeatletter
\usepackage[toc,page,header]{appendix}
\usepackage{minitoc}
\usepackage{ifthen}
\newboolean{doublecolumn}
\setboolean{doublecolumn}{false}
\newboolean{arxiv}
\setboolean{arxiv}{true}
\usepackage{pgfplots}
\usepackage{pdflscape}
\usepackage{hyperref}
\usepackage{cleveref}
\usepackage{authblk}
\pgfplotsset{compat=1.15}

\newcommand{\assign}{:=}

\newcommand{\cdummy}{\cdot}
\newcommand{\tmtextbf}[1]{\text{{\bfseries{#1}}}}
\newcommand{\tmop}[1]{\ensuremath{\operatorname{#1}}}

\newcommand{\minf}[1]{\underset{#1}{\text{minimize}}}
\newcommand{\ogm}{{\texttt{OGD}}}
\newcommand{\osgm}{{\texttt{OSGM}}}
\newcommand{\osgmrx}{\texttt{OSGM-R}}
\newcommand{\osgmrzx}{\texttt{OSGM-RZ}}
\newcommand{\osgmgx}{\texttt{OSGM-G}}
\newcommand{\osgmhx}{\texttt{OSGM-H}}
\newcommand{\gd}{{\texttt{GD}}}
\newcommand{\optdgd}{{\texttt{OptDiagGD}}}
\newcommand{\adagrad}{{\texttt{AdaGrad}}}
\newcommand{\sagd}{{\texttt{SAGD}}}
\newcommand{\agd}{{\texttt{AGD}}}
\newcommand{\mathd}{\mathrm{d}}
\global\long\def\vertiii#1{\left\vert \kern-0.25ex  \left\vert \kern-0.25ex  \left\vert #1\right\vert \kern-0.25ex  \right\vert \kern-0.25ex  \right\vert }%

\global\long\def\diam{\mathrm{diam}}%

\global\long\def\argmin{\operatornamewithlimits{arg\,min}}%

\global\long\def\and{\mathrm{and}}%

\global\long\def\Rbb{\mathbb{R}}%

\global\long\def\Acal{\mathcal{A}}%

\global\long\def\Lcal{\mathcal{L}}%

\global\long\def\Mcal{\mathcal{M}}%

\global\long\def\Ncal{\mathcal{N}}%

\global\long\def\Ocal{\mathcal{O}}%

\global\long\def\Pcal{\mathcal{P}}%

\global\long\def\Xcal{\mathcal{X}}%

\theoremstyle{plain}
\newtheorem{lem}{\protect\lemmaname}[section]
\theoremstyle{remark}
\newtheorem{rem}{\protect\remarkname}
\theoremstyle{plain}
\newtheorem{thm}{\protect\theoremname}[section]
\theoremstyle{plain}
\newtheorem{prop}{\protect\propositionname}[section]
\providecommand{\corollaryname}{Corollary}
\theoremstyle{plain}
\newtheorem{coro}{\protect\corollaryname}[section]
\theoremstyle{plain}

\theoremstyle{plain}
\newtheorem{definition}{\protect\definitionname}[section]

\providecommand{\lemmaname}{Lemma}
\providecommand{\remarkname}{Remark}
\providecommand{\theoremname}{Theorem}
\providecommand{\examplename}{Example}
\providecommand{\propositionname}{Proposition}
\providecommand{\definitionname}{Definition}

\crefdefaultlabelformat{#2\textbf{#1}#3} % <-- Only #1 in \textbf
\crefname{section}{\textbf{section}}{\textbf{sections}}
\Crefname{section}{\textbf{Section}}{\textbf{Sections}}
\crefname{thm}{\textbf{theorem}}{\textbf{theorems}}
\Crefname{thm}{\textbf{Theorem}}{\textbf{Theorems}}
\crefname{lem}{\textbf{lemma}}{\textbf{lemmas}}
\Crefname{lem}{\textbf{Lemma}}{\textbf{Lemmas}}
\crefname{prop}{\textbf{proposition}}{\textbf{propositions}}
\Crefname{prop}{\textbf{Proposition}}{\textbf{Propositions}}
\crefname{algorithm}{\textbf{algorithm}}{\textbf{algorithms}}
\Crefname{algorithm}{\textbf{Algorithm}}{\textbf{Algorithms}}
\crefname{coro}{\textbf{Corollary}}{\textbf{corollaries}}
\Crefname{coro}{\textbf{Corollary}}{\textbf{corollaries}}
\crefname{definition}{\textbf{Definition}}{\textbf{definitions}}
\Crefname{definition}{\textbf{Definition}}{\textbf{definitions}}
\crefname{table}{\textbf{Table}}{\textbf{tables}}
\Crefname{table}{\textbf{Table}}{\textbf{tables}}
\crefname{figure}{\textbf{Figure}}{\textbf{figures}}
\Crefname{figure}{\textbf{Figure}}{\textbf{figures}}

\newcommand{\rk}{r_{k}}
\newcommand{\rxz}{r^z_{x}}
\newcommand{\rxkz}{r^z_{x^k}}

\newcommand{\YC}[1]{ }
\renewcommand{\YC}[1]{\textcolor{blue}{[YC: #1]}}

\begin{document}

\title{Gradient Methods with Online Scaling}

\author[1]{Wenzhi Gao\thanks{gwz@stanford.edu}}
\author[2]{Ya-Chi Chu\thanks{ycchu97@stanford.edu}}
\author[1,3]{Yinyu Ye\thanks{yyye@stanford.edu}}
\author[1,3]{Madeleine Udell\thanks{udell@stanford.edu}}
\affil[1]{ICME, Stanford University}
\affil[2]{Department of Mathematics, Stanford University}
\affil[3]{Department of Management Science and Engineering, Stanford University}

\maketitle

 \begin{abstract}
We introduce a framework to accelerate the convergence of gradient-based methods with online learning. 
The framework learns to scale the gradient at each iteration through an online learning algorithm and provably accelerates gradient-based methods asymptotically. 
In contrast with previous literature, where convergence is established based on worst-case analysis, our framework provides a strong convergence guarantee with respect to the optimal scaling matrix for the \textit{iteration trajectory}.
For smooth strongly convex optimization, our results provide an $\Ocal(\kappa^\star \log(1/\varepsilon)$) complexity result, where $\kappa^\star$ is the condition number achievable by the optimal preconditioner, improving on the previous $\Ocal(\sqrt{n}\kappa^\star \log(1/\varepsilon))$ result. In particular, a variant of our method achieves superlinear convergence on convex quadratics. For smooth convex optimization, we show for the first time that the widely-used hypergradient descent heuristic improves on the convergence of gradient descent.
\end{abstract}
\section{Introduction} \label{sec:intro}

We consider the unconstrained smooth strongly convex optimization problem
\begin{eqnarray*}
  \minf{x \in \mathbb{R}^n} & f (x), & 
\end{eqnarray*}
where $f (x) : \mathbb{R}^n \rightarrow \mathbb{R}$ is $L$-smooth and
$\mu$-strongly convex with $f (x^{\star}) \assign \min_x f (x) > - \infty$.
 It is well-known that gradient descent with stepsize $1 / L$ converges with
iteration complexity $\mathcal{O} (\kappa \log (1 / \varepsilon))$, where
$\kappa = L / \mu$ is the condition number of the problem. Two major
techniques have been developed in the literature to accelerate gradient descent. One is
to improve the dependence on $\kappa$ through Nesterov's fast gradient method
{\cite{necoara2019linear,nesterov2013introductory}};
the other is through preconditioning: a positive definite matrix stepsize $P$, known as preconditioner, premultiplies the gradient to improve convergence:
\[ x^{k + 1} = x^k - P \nabla f (x^k). \]
Preconditioning has been a standard tool in convex
optimization and numerical linear algebra to improve convergence of gradient descent {\cite{li2017preconditioned,maddison2021dual,li2016preconditioned,frangella2022sketchysgd,frangella2023promise}} or other iterative methods {\cite{saad2003iterative}}, and
it is closely related to the well-known adaptive gradient methods
{\cite{duchi2011adaptive,kingma2014adam}}, either for online learning or for a
general optimization problem. Some recent results quantify the effect of
adaptive methods on problem conditioning {\cite{das2024towards}}. In the
context of machine learning, adaptively choosing a preconditioner is also
relevant to hyperparameter tuning
{\cite{hospedales2021meta,necoara2019linear}}, especially choosing a learning rate schedule {\cite{defazio2024road}}.\\

Despite the great empirical success of adaptive methods in practice, they usually cannot improve the theoretical complexity as a function of the condition number. 
Recently, {\cite{kunstner2024searching}} showed that hypergradient, the gradient of the optimization objective with respect to the preconditioner, can be used to improve the convergence of gradient descent. The idea in {\cite{kunstner2024searching}}  is to use a cutting plane subroutine to update the (diagonal) preconditioner, and an $\mathcal{O} ( \sqrt{n}
\kappa^{\star} \log (1 / \varepsilon) )$ complexity result is obtained, where $n$ is the variable dimension and $\kappa^\star$ is the condition number of the optimally preconditioned problem. Although the result in {\cite{kunstner2024searching}} is
dimension-dependent and requires a nontrivial subroutine to update the
preconditioner, it provides a valuable direction to improve the
performance of first-order adaptive methods theoretically. Whether a simple adaptive first-order method can achieve $\mathcal{O} (\kappa^{\star} \log (1 / \varepsilon))$ complexity or even stronger guarantees
 remains open.\\

This paper answers this question affirmatively by proposing the online scaled gradient method, a framework that accelerates gradient-based methods through online convex optimization. 

\paragraph{Contributions.}
\begin{itemize}[leftmargin=10pt]
  \item We develop a framework that accelerates gradient-based algorithms through online learning. Unlike previous work, our framework guarantees convergence with respect to the scaling matrix optimized for the iteration trajectory, rather than the worst-case analysis. 
  
  \item We propose a simple adaptive first-order gradient-based method with asymptotic $\mathcal{O} (\kappa^{\star} \log (1 / \varepsilon))$ complexity, where $\kappa^\star$ is the optimal condition number achievable by the optimal preconditioner, improving on the $\mathcal{O} (\sqrt{n}\kappa^{\star} \log (1 / \varepsilon))$ complexity in the previous literature. In particular, one realization of our framework achieves superlinear convergence on strongly convex quadratics using first-order information.
  
  \item For the first time, we prove that the hypergradient heuristic improves the convergence of gradient descent.
  
\end{itemize}

\subsection{Related literature}

\paragraph{Preconditioned iterative methods.} Preconditioning is a well-established technique to enhance the convergence of iterative algorithms in both optimization \cite{frangella2023promise, o2016conic, applegate2021practical, deng2024enhanced} and numerical linear algebra \cite{saad2003iterative, qu2024optimal, gao2023scalable, frangella2023randomized, doan2011numerical}. By applying a linear transformation to the optimization variables, preconditioning aims to reduce the heterogeneity of the optimization landscape. Recent research has focused on understanding the properties of optimal preconditioners \cite{qu2024optimal, gao2023scalable, jambulapati2020fast}. While these methods demonstrate empirical success, identifying a good preconditioner can be computationally intensive and often depends on the specific structure of the problem. 

\paragraph{Hypergradient descent heuristic.} Our method is closely related to the hypergradient descent heuristic   \cite{gunes2018online,maclaurin2015gradient,chandra2022gradient}, which updates the stepsize (hyperparameters) using the gradient of the optimization objective with respect to it. 
Despite strong empirical results \cite{gunes2018online, chandra2022gradient}, the theoretical understanding of hypergradient descent remains limited. The existing results \cite{rubio2017convergence} are unable to fully justify the observed improvements \cite{kunstner2024searching}. Recently, \cite{kunstner2024searching} introduced a novel multi-dimensional backtracking approach that uses hypergradients to generate separating hyperplanes in the space of candidate preconditioners. This work provides the first theoretical justification for hypergradient descent. However, it solves a small cutting plane subproblem at every step and incurs a $\sqrt{n}$ dependence on the problem dimension. No theoretical proof exists that the original hypergradient descent heuristics accelerate gradient-based methods. And our paper provides the first proof that quantifies the acceleration effect of the hypergradient descent heuristic.

\paragraph{Adaptive first-order methods.} Adaptive stepsize is a well-established technique to enhance the convergence of optimization algorithms. The most notable example is \texttt{AdaGrad} \cite{duchi2011adaptive, mcmahan2010adaptive}, which provides strong theoretical guarantees in the context of online convex optimization. Other methods, such as \texttt{Adam} \cite{kingma2014adam} and \texttt{RMSProp} \cite{hinton2012neural}, have demonstrated competitive empirical performance, though they generally yield weaker online regret bounds. Our approach also leverages online learning techniques to accelerate gradient-based methods, with a particular focus on improving the dependence on problem conditioning. We show that online learning applied to the hypergradient gives optimal guarantees on problem conditioning.

\paragraph{Learning to optimize and meta-learning.} The learning to optimize \cite{li2016learning,chen2022learning} and meta-learning \cite{hospedales2021meta, chen2024online, finn2019online} literatures also use online learning to improve algorithm performance. These approaches are typically designed to solve a sequence of related optimization problems, providing performance guarantees across multiple tasks. In contrast, our work applies online learning to improve first-order methods over the course of solving a single optimization instance.
\section{Background and preliminaries}

\paragraph{Notations.}Throughout the paper, we use $\| \cdummy \|$ to denote vector Euclidean norm or matrix spectral norm, and $\langle \cdummy, \cdummy \rangle$ to denote Euclidean
inner product. Letters $A$ and $a$ denote matrices and vectors, respectively.
$\| A \|_F = \sqrt{\sum_{i j} a_{i j}^2}$ denotes the matrix Frobenius norm. Given two vectors $a, b$ of the same dimension, $a \odot b$ denotes their element-wise Hadamard product.
The Clarke subdifferential of a function $f (x)$ at $x$ is defined by $\partial f (x) \assign \{ v \in \mathbb{R}^n : f (y) \geq f (x) + \langle v, y - x \rangle + o(\|x-y\|), y \rightarrow x\}$. 
We use $f' (x) \in \partial f (x)$ to denote a subgradient. If $f (x)$ is differentiable at $x$, $\partial f (x) = \{ \nabla f (x) \}$. 
For symmetric matrices $A, B$, $A \succeq B$ if $A - B \in
\mathbb{S}^n_+$ is positive semidefinite. The condition number for an $L$-smooth
and $\mu$-strongly convex problem is defined by $\kappa = L / \mu$. Given a closed
convex set $C$, $\Pi_C [x]$ denotes the orthogonal projection of $x$
onto $C$. We use $\Lcal_\alpha := \{x: f(x) \leq \alpha\}$ to denote the $\alpha$-sublevel set of $f$ and $\Xcal^\star$ to denote the optimal set of $f$. $x$ is an $\varepsilon$-optimal solution if $f(x) \leq f(x^\star) + \varepsilon$ and $x$ is an $\varepsilon$-critical point if $\|\nabla f(x)\| \leq \varepsilon$.

\subsection{Assumptions}

We make the following two assumptions throughout the paper.

\begin{enumerate}[leftmargin=35pt,itemsep=2pt,label=\textbf{A\arabic*:},ref=\rm{\textbf{A\arabic*}},start=1]
\item $f (x)$ is $L$-smooth. $| f (x) - f (y) - \langle \nabla f (y),
x - y \rangle | \leq \tfrac{L}{2} \| x - y \|^2$ \label{A1}

\item $f (x)$ is $\mu$-strongly convex. $f (x) - f (y) - \langle
\nabla f (y), x - y \rangle \geq \frac{\mu}{2} \|x - y\|^2$ \label{A2}
\end{enumerate}

We also assume that $f$ is twice-differentiable for simplicity. However, our algorithm does not necessarily require twice-differentiability
to work. In addition, $\mu$-strong convexity can be relaxed to weaker
conditions such as convexity with quadratic growth {\cite{necoara2019linear}}.

\subsection{Preconditioned and scaled gradient method}

It is well-known that under \ref{A1} and \ref{A2}, the vanilla gradient descent
\[ x^+ = x - \tfrac{1}{L} \nabla f (x) \]
achieves $\mathcal{O} (\kappa \log (1 / \varepsilon))$ complexity
{\cite{garrigos2023handbook}}. 
The dependence on the condition number $\kappa$ is unfortunate since the condition number can be very large and substantially slows down convergence.
Two techniques are often used to improve
dependence on $\kappa$. One is through Nesterov's fast gradient method
{\cite{necoara2019linear}}, which achieves $\mathcal{O} ( \sqrt{\kappa} \log (1 /
\varepsilon) )$ complexity; the other is through
preconditioning
{\cite{frangella2022sketchysgd,li2017preconditioned,frangella2023promise}},
which replaces the scalar stepsize by some positive semidefinite matrix $P \in
\mathbb{S}^n_{+}$
\begin{equation} \label{eqn:pgd-update}
	x^+ = x - P \nabla f (x).
\end{equation}
Typically, $P$ is chosen to be diagonal and positive definite. 
This paper also considers updates of the form \eqref{eqn:pgd-update} but allows $P$ to be an arbitrary matrix of proper dimension from some closed convex set
$\mathcal{P}$. Moreover, $P$ is allowed to vary across
iterations: 
\begin{equation} \label{eqn:sgm-update}
	x^{k+1} = x^k - P_k \nabla f (x^k) 
\end{equation}

To differentiate our method from standard preconditioning techniques, which generally consider only positive definite preconditioners, we call our method the
scaled gradient method: $P$ serves as a (not necessarily positive definite or symmetric) scaling matrix. Preconditioned gradient
descent can be viewed as a special case of the scaled gradient method. We define
$\mathcal{P}_+ :=\mathcal{P} \cap \mathbb{S}^n_+$ when we switch to the context of preconditioned gradient descent, and without loss of generality we assume $0 \in \Pcal, L^{-1}I\in \Pcal$ and that $\Pcal$ is bounded.

\begin{enumerate}[leftmargin=35pt,itemsep=2pt,label=\textbf{A\arabic*:},ref=\rm{\textbf{A\arabic*}},start=3]
   \item Closed convex set $\Pcal$ satisfies $0\in \Pcal$, $L^{-1}I\in \Pcal$ and $\diam(\Pcal) \leq D$. \label{A3}
\end{enumerate}

\subsection{Monotone descent oracle} \label{sec:monotone-oracle}
When $P$ is not positive definite, a scaled gradient update will not necessarily decrease the function value. To guarantee convergence under weak assumptions, 
the scaled gradient method optionally uses a \emph{monotone descent
oracle} $\mathcal{M}$, defined below.

\begin{definition} \label{def:monotone-oracle}
  Given the scaled gradient update $x^+ = x - P \nabla f (x)$, $\Mcal_{\varphi, P}: \Rbb^n \rightarrow \Rbb^n$ is called a monotone descent oracle associated with the scaled gradient update and measure $\varphi$ if its output $\mathcal{M}_{\varphi, P} (x)$ satisfies
  \[ \varphi (\mathcal{M}_{\varphi, P} (x)) \leq \min \{ \varphi (x), \varphi (x^+) \} . \]
\end{definition}

We use $\Mcal(x)$ to denote the oracle when the context is clear. Three typical realizations of $\mathcal{M}$ are as follows.
\begin{itemize}[leftmargin=10pt]
  \item \emph{Line-search.}
  
  $\mathcal{M} (x) = x + \alpha (x^{+} - x)$ such that $\varphi (\mathcal{M} (x))
  \leq \varphi (x)$. Additional regularity conditions, such as $\mathcal{P} = \Pcal_+$, are required to ensure that line-search stops in a finite number of steps.
  
  \item \emph{Steepest descent.}
  
  $\mathcal{M} (x) = x + \alpha (x^{+} - x)$ and $\alpha = \argmin_\alpha \varphi(x + \alpha (x^{+} - x))$. It applies to simple functions such as quadratic.
  
  \item \emph{Simple comparison.}
  
  $\mathcal{M} (x) = x^{+}$ if $\varphi (x^{+}) \leq \varphi (x)$. Otherwise
  $\mathcal{M} (x) = x$. It takes one extra measure evaluation.
\end{itemize}

\section{Online Scaled Gradient Methods} \label{sec:online}

This section introduces our main methodology, which relates the scaled gradient method to online convex optimization with $P$ as the decision variable.

\subsection{Scaled gradient method and online learning} \label{sec:sgm-ol}

Let $\varphi (x)$ be a non-negative measure or potential energy that characterizes the optimality of $x$. For example, function value gap $\varphi (x) = f (x) - f (x^{\star})$ and gradient norm $\varphi (x) = \| \nabla f
(x) \|$ are common measures. 
The progress of an algorithm at step $K+1$ with respect to measure $\varphi$ can be expressed as the telescoping product
\begin{equation} \label{eqn:tele-prod} 
   \varphi (x^{K + 1}) = \varphi(x^1) \prod_{k=1}^K \frac{\varphi(x^{k+1})}{\varphi(x^{k})}.
\end{equation}
Then the arithmetic-geometric mean inequality upper-bounds $\varphi(x^{K+1})$:

\begin{thm}\label{thm:online-to-conv-ratio}
Given a non-negative function $\varphi (x) : \mathbb{R}^n \rightarrow
\mathbb{R}_+$ and a sequence of iterations $\{ x^k \}$,
\begin{equation*}
   \varphi (x^{K + 1}) \leq \varphi (x^1) \big( \tfrac{1}{K} \textstyle \sum_{k = 1}^K \frac{\varphi (x^{k + 1})}{\varphi (x^k)} \big)^K. 
\end{equation*}
\end{thm}

The quantity $\frac{1}{K} \textstyle \sum_{k = 1}^K \frac{\varphi (x^{k + 1})}{\varphi (x^k)}$ on the right-hand side is the averaged contraction factor across all previous iterates: a smaller contraction factor ensures stronger convergence. Suppose the iterates $\{ x^k \}_{k \geq 2}$ are generated by the scaled gradient method in \eqref{eqn:sgm-update}. Then
\begin{equation} \label{eqn:scaled-grad-bound}
   \tfrac{1}{K} \textstyle \sum_{k = 1}^K \frac{\varphi (x^{k + 1})}{\varphi (x^k)} = \frac{1}{K}  \sum_{k = 1}^K \frac{\varphi (x^k - P_k \nabla f (x^k))}{\varphi (x^k)}.
\end{equation}
To maximize the progress in the scaled gradient method, we aim to minimize the quantity in \eqref{eqn:scaled-grad-bound} over the choice of scaling matrices $P_k$ with online learning.  
We will show that online convex optimization can learn a sequence of $\{P_k\}$ that asymptotically accelerates gradient-based methods. Define the \emph{surrogate loss} $$\ell_x (P) \assign \tfrac{\varphi (x - P \nabla f (x))}{\varphi(x)}$$ with respect to the measure $\varphi$. Note that $\ell_{x^k}$ only depends on $x^1$ and all previous scaling matrices $\{P_j\}_{j \leq k - 1}$.

Online learning generates a sequence $\{P_k\}$ such that the cumulative regret is bounded by $\rho_K$: 
\begin{equation} \label{eqn:def-regret}
\textstyle \sum_{k = 1}^K \ell_{x^k} (P_k) - \displaystyle \min_{P \in \mathcal{P}} \textstyle \sum_{k = 1}^K \ell_{x^k} (P) \leq \rho_K.
\end{equation}
Existing results in online optimization can guarantee sublinear regret if the losses $\{\ell_{x^k}\}$ are convex and are either Lipschitz continuous or have Lipschitz continuous gradient \cite{orabona2019modern}. 
In this case, we say the family of surrogate losses $\{\ell_{x^k}\}$ is \emph{online-learnable}.
The definition of regret $\rho_K$ and \Cref{thm:online-to-conv-ratio} imply
\begin{equation} \label{eqn:progress-bound-by-regret}
   \varphi (x^{K + 1}) \leq \varphi (x^1) \big( 
   \tfrac{1}{K} \textstyle \sum_{k = 1}^K \ell_{x^k} (P)  \big)^K \leq \varphi (x^1) \big( \displaystyle \min_{P \in \mathcal{P}} \tfrac{1}{K} \textstyle \sum_{k = 1}^K \ell_{x^k} (P) + \tfrac{\rho_K}{K} \big)^K .
\end{equation}
When the regret $\rho_K$ grows sublinearly in $K$, the bound in \eqref{eqn:progress-bound-by-regret} suggests that for large enough $K$,
\begin{equation*}
   \varphi (x^{K + 1}) \leq \varphi (x^1)
   \big( \displaystyle \min_{P \in \mathcal{P}} \tfrac{1}{K} \textstyle \sum_{k = 1}^K \ell_{x^k} (P) + \tfrac{\rho_K}{K} \big)^K 
   \approx \varphi (x^1) \big( \displaystyle \min_{P \in \mathcal{P}} \tfrac{1}{K} \textstyle \sum_{k = 1}^K \ell_{x^k} (P) \big)^K.
\end{equation*}

This result is powerful: it suggests that a scaled gradient method, in the long run, can achieve convergence that is competitive with any fixed scaling matrix optimized for the \emph{iteration trajectory}. To the best of our knowledge, this trajectory-based convergence guarantee is rare in the literature. Moreover, as long as there exists some pre-specified scaling matrix $P^\star$ (or simply stepsize $P^\star = \alpha I$) such that $\ell_{x}(P^\star)\leq \theta^\star < 1$ for any $x$, we obtain the global convergence guarantee
\begin{align}
   \varphi (x^{K + 1}) \leq{} \varphi(x^1) (\theta^\star + \tfrac{\rho_K}{K})^K \approx \varphi(x^1) (\theta^\star)^K
\end{align}
The algorithm, which updates the scaling matrix $P_k$ on the fly, is called a realization of the \emph{online scaled gradient method} (\osgm).

\subsection{Framework of online scaled gradient method}

The online scaled gradient method is determined by the components below:
\begin{itemize}[leftmargin=15pt]
  \item \emph{Optimality measure}. A measure $\varphi$ to characterize the convergence of {\osgm}.
  
  \item \emph{Surrogate loss}. A surrogate loss $\ell_x(P)$ that relates  $\varphi$ with an online learning problem in $P$. 
  
  \item \emph{Online learning algorithm}. An online learning algorithm $\mathcal{A}$ that guarantees sublinear regret for $\sum_k \ell_{x^k}(P_k)$.
  
  \item (Optional) \emph{Monotone oracle}. An oracle $\Mcal$ (\Cref{def:monotone-oracle}) that guarantees monotonicity.
  
  \item (Theoretical purpose) \emph{Hindsight scaling matrix.} A hindsight scaling matrix $P^\star$ to ensure global convergence.
  
\end{itemize}
The tuple $(\varphi, \ell, \mathcal{A}, \Mcal)$ determines a realization
of the online scaled gradient method (\osgm, \Cref{alg:osgm}).

\begin{algorithm}[h]
{\textbf{input $x^1, P_1, \varphi, \ell, \mathcal{A}, \Mcal$}}\\
\For{k =\rm{ 1, 2,...}}{
\eIf{$\Mcal = \varnothing $}{
$x^{k + 1} = x^{k} - P_k \nabla f(x^k)$
}{
$x^{k + 1} = \Mcal_{\varphi, P_k}(x^{k})$
}
$P_{k + 1} = \Acal(\ell_{x^k}, P_k) $
}
{\textbf{output} $x^\text{best}$ with minimum objective value}
\caption{Online scaled gradient method (\osgm) \label{alg:osgm}}
\end{algorithm} 

The surrogate loss $\ell_x$ is the key component of this framework.
Finding a good surrogate, for given assumptions on the function to be optimized, is the challenging aspect of the theoretical analysis. In the rest of the paper, we provide several realizations of the framework for different function classes, summarized in \Cref{tab:summary}.

\begin{table}[h]
\centering
\renewcommand{\arraystretch}{1}
  \begin{tabular}{ccccccc}
    \toprule
    $\varphi (x)$ & Surrogate $\ell$ & Strong convexity & $\mathcal{A}$ &
    $\mathcal{M}$ oracle & Complexity & Reference\\
    \midrule
    $f (x) - f (x^{\star})$ & $r_x (P) = \frac{f (x^+) - f
    (x^{\star})}{f (x) - f (x^{\star})}$ & Yes & & Optional & $\mathcal{O}
    ( \kappa^{\star} \log ( \tfrac{1}{\varepsilon} )
    )$ & \Cref{sec:ratio-surrogate}\\
%    \midrule
    $\| \nabla f (x) \|$ & $g_x (P) = \frac{\| \nabla f (x^+) \|}{\|
    \nabla f (x) \|}$ & Yes &  & Required & $\mathcal{O} ( \lambda^{\star}
    \log ( \tfrac{1}{\varepsilon} ) )$ & \Cref{sec:gnorm-surrogate}\\
%    \midrule
    \multirow{2}{*}{$f (x) - f (x^{\star})$} & \multirow{2}{*}{$h_x (P) = \frac{f (x^+) - f (x)}{\|
    \nabla f (x) \|^2}$} & {Yes} & \multirow{-2.5}{*}{\ogm} & Required & $\mathcal{O} (
    \frac{1}{2 \mu\gamma^{\star}} \log ( \tfrac{1}{\varepsilon} ) )$ & \multirow{2}{*}{\Cref{sec:hypergrad-surrogate}}\\
    & &No & & Required & $\mathcal{O} (
     \tfrac{1}{\gamma^{\star}\varepsilon} )$ & \\
    \bottomrule
  \end{tabular}
  \caption{Realizations of \osgm. \ogm:
  online (sub)gradient method.
  $\kappa^{\star}, \lambda^{\star}, \gamma^{\star}$ and their optimal scaling matrix
  $P^{\star}$ will be defined in the next sections.}
  \label{tab:summary}
\end{table}
\section{Function value ratio surrogate}\label{sec:ratio-surrogate}

The first surrogate, function value {\textit{ratio surrogate}}, is defined as follows:
\begin{equation}\label{eqn:ratio-surrogate-loss}
	r_x (P) \assign \tfrac{f (x^+) - f (x^{\star})}{f (x) - f (x^{\star})} =
   \tfrac{f (x - P \nabla f (x)) - f (x^{\star})}{f (x) - f (x^{\star})} . 
\end{equation}
The ratio surrogate $r_x$ measures the contraction factor of the function value gap between two consecutive \osgm~steps. We assume strong convexity ($\mu > 0$) throughout this section.
Without loss of generality, we assume all the intermediate iterates generated by the algorithm satisfy $f(x) > f(x^\star)$ so that $r_x$ is well-defined; otherwise, $x \in \Xcal^\star$ and we can immediately stop the algorithm.  The ratio surrogate $r_x$ assumes the optimal value $f (x^{\star})$ is known, and this assumption will be relaxed later in this section. The monotone oracle is optional for $r_x$. We present the results without a monotone oracle: $\Mcal = \varnothing$.

\subsection{Surrogate loss}
The function value ratio $r_x$ in \eqref{eqn:ratio-surrogate-loss} can be viewed as a surrogate loss since its average along \osgm~iterates serves as an upper bound on the function value gap. By substituting $\varphi (x) = f (x) - f (x^{\star})$ in \Cref{thm:online-to-conv-ratio} and plugging in the definition of $r_x$, \osgm~iterates $\{x^k\}$ satisfy the following bound:

\begin{lem}[Surrogate loss and measure] \label{lem:rx-surrogate-measure}
For all $K \geq 1$, the online scaled gradient method satisfies
\begin{equation} \label{eqn:rx-surrogate-measure}
f (x^{K + 1}) - f (x^{\star}) \leq ( f (x^1) - f (x^{\star}) ) \big( \tfrac{1}{K} \textstyle \sum_{k = 1}^K r_{x^k} (P_k) \big)^K. 
\end{equation}
\end{lem}
The ratio surrogate $r_x$ inherits several important properties from $f$, which we summarize in \Cref{prop:rx-learnability}. These properties are crucial for establishing the online learnability of $\{r_{x^k}\}$ later in \Cref{sec:ratio-ol-alg}.

\begin{prop}[Properties of $r_x$] \label{prop:rx-learnability}
Under \ref{A1} and \ref{A2}, for any fixed $x \not \in \Xcal^\star $, the surrogate loss $r_x (P)$ defined in \eqref{eqn:ratio-surrogate-loss} is convex, non-negative, and $2 L^2$-smooth as a function in $P$. In addition, the derivative of $r_x$ takes the form
\begin{equation} \label{eqn:grad-ratio}
   \nabla r_{x} (P) = - \tfrac{\nabla f (x - P \nabla f(x)) \nabla f (x)^{\top}}{f (x) - f (x^{\star})}.
\end{equation}
\end{prop}

\subsection{Online learning algorithm} \label{sec:ratio-ol-alg}
Online gradient descent is known to ensure sublinear regret for a family of smooth, convex, and lower-bounded losses \cite{orabona2019modern}, which is the case for ratio surrogate loss $r_x$ by \Cref{prop:rx-learnability}.
We tailor the classical $L^\star$ regret bound from online convex optimization literature \cite{orabona2019modern} to our settings in \Cref{lem:rx-learnability} below.

\begin{lem}[Learnability] \label{lem:rx-learnability}
Given \ref{A1}, \ref{A2}, and the ratio surrogate losses $\{r_{x^k}\}$,  online gradient descent
\begin{equation}\label{eqn:rx-ogd}
	P_{k + 1} = \Pi_{\mathcal{P}} [P_k - \eta \nabla r_{x^k} (P_k)]
\end{equation}
with stepsize $\eta \leq 1/(4 L^2)$ generates a sequence of scaling matrices $\{P_k\}_{k \geq 2}$ such that 
\begin{equation} \label{eqn:rx-regret-bound}
   \textstyle \sum_{k = 1}^K r_{x^k} (P_k) - \textstyle \sum_{k = 1}^K r_{x^k} (P) \leq \tfrac{1}{\eta} \| P - P_1 \|^2_F + 4 L^2 \eta \textstyle \sum_{k = 1}^K r_{x^k} (P) \quad \text{for any } P \in \mathcal{P}.
\end{equation}
In particular, if \ref{A3} is further assumed, the choice of stepsize $\eta = \min\big\{ \tfrac{1}{4 L^2}, \tfrac{D}{2L (1 + LD) \sqrt{K}} \big\}$ ensures
\begin{equation} \label{eqn:r-regret-bound-D}
   \textstyle \sum_{k = 1}^K r_{x^k} (P_k) - \displaystyle \min_{P \in \mathcal{P}} \textstyle \sum_{k = 1}^K r_{x^k} (P) 
   \leq  \rho_K := \max\big\{ 4 LD (1 + LD) \sqrt{K}, 8 L^2 D^2 \big\}.
\end{equation}
\end{lem}

\begin{rem}
The relation \eqref{eqn:rx-ogd} suggests additional complexity from a rank-one update with an orthogonal projection. But as we will discuss in \Cref{sec:practical}, we can choose $\Pcal$ to have arbitrary sparsity (e.g., diagonal), and it is only necessary to update the nonzero elements. Moreover, the orthogonal projection is often easy to compute since we do not require $P_k$ to be positive semidefinite or symmetric.
\end{rem}

\subsection{Algorithm design and analysis}

We now state a realization of {\osgm} with the ratio surrogate loss $r_{x}$, denoted by \osgmrx. 
We choose the optimality measure $\varphi$, the surrogate loss $\ell$, and the online learning algorithm $\mathcal{A}$ to be
\begin{equation*}
    \varphi(x) \assign f(x) - f(x^\star), \quad \ell_{x} (P) \assign r_{x} (P), \quad \mathcal{A} \assign \text{online gradient descent in \eqref{eqn:rx-ogd}},
\end{equation*}
and the monotone oracle $\Mcal$ is optional.
\Cref{alg:osgm-rx} presents \osgmrx~without the monotone oracle.

\begin{algorithm}[h]
{\textbf{input} $x^1, P_1 \in \Pcal, \text{online gradient stepsize }\eta > 0$}\\
\For{k =\rm{ 1, 2,...}}{
$x^{k + 1} = x^{k} - P_k \nabla f(x^k)$\\
$P_{k+1} = \Pi_{\Pcal} [P_k - \eta \nabla r_{x^k} (P_k) ]$

}
{\textbf{output} $x^\text{best}$ with minimum objective value}
\caption{Online scaled gradient method with ratio surrogate (\osgmrx) \label{alg:osgm-rx}}
\end{algorithm} 

Combining \Cref{lem:rx-surrogate-measure} and \Cref{lem:rx-learnability}, \Cref{thm:rx-trajectory-opt} characterizes the trajectory-based convergence behavior of {\osgmrx}.

\begin{thm}[Trajectory-based convergence] \label{thm:rx-trajectory-opt}
Under \ref{A1} to \ref{A3}, \Cref{alg:osgm-rx} (\osgmrx) with \\$\eta = \min\big\{ \tfrac{1}{4 L^2}, \tfrac{D}{2L (1 + LD) \sqrt{K}} \big\}$ satisfies
\begin{equation} \label{eqn:rx-trajectory-opt}
   f (x^{K + 1}) - f (x^{\star}) 
   \leq (f (x^1) - f (x^{\star})) ( \theta^{\star}_K +
   \tfrac{\rho_K}{K} )^K,
\end{equation}
where $\theta^{\star}_K \assign  \min_{P \in \mathcal{P}} \tfrac{1}{K} \textstyle \sum_{k = 1}^K r_{x^k} (P)$ and $\rho_K = \max\big\{ 4 LD (1 + LD) \sqrt{K}, 8 L^2 D^2 \big\}$ is defined in \eqref{eqn:r-regret-bound-D}.
\end{thm}

From \eqref{eqn:rx-trajectory-opt}, when $K$ is large enough, $\tfrac{\rho_K}{K}$ is negligible, and \osgmrx~behaves like an algorithm with linear convergence rate $\theta^{\star}_K$.  Note that $\theta^{\star}_K$ is based on the optimization trajectory, and the behavior of {\osgmrx} is competitive with the scaling matrix that minimizes $\theta^{\star}_K$. To our knowledge, this trajectory-based convergence guarantee is new in the literature. To show global convergence, we need to show that $\theta^{\star}_K < 1$, and this fact follows from the existence of the \emph{optimal preconditioner}.

\subsubsection{Hindsight and global convergence}
Define $P_r^{\star}$ to be the ratio scaling matrix that solves the following semidefinite optimization problem
\begin{equation}\label{eqn:optprecond}
	\kappa^{\star} \assign \min_{P \in \mathcal{P}_+} \kappa \quad
   \text{subject to} \quad \tfrac{1}{\kappa} I \preceq P^{1 / 2} \nabla^2 f (x)
   P^{1 / 2} \preceq I \quad \text{for all~~} x.
\end{equation}
$P_r^{\star}$ is known in the literature {\cite{kunstner2024searching, qu2024optimal}} as the \textit{universal optimal preconditioner}.
The optimal value $\kappa^{\star}$ is called the \emph{optimal condition number} with respect to subset $\mathcal{P}_+ = \Pcal \cap \mathbb{S}^n_{+} $. Since \ref{A3} assumes $L^{-1} I\in \Pcal$, we have $\kappa^\star \leq \kappa$ and preconditioned gradient descent with preconditioner $P_r^{\star}$ converges as if the condition number of the underlying minimization problem is reduced from $\kappa = L / \mu$ to $\kappa^{\star}$. A standard argument using the descent lemma and strong convexity of $f(x)$ (\ref{A2}) ensures
\begin{equation} \label{eqn:opt-contraction-factor}
f(x - P_r^\star \nabla f(x)) - f(x^\star) \leq (1 - \tfrac{1}{\kappa^\star}) (f(x) - f(x^\star)) \quad \text{for all~~} x,
\end{equation}
which can be equivalently expressed in terms of the ratio surrogate loss $r_x$ in the lemma below:

\begin{lem}[Hindsight] \label{lem:rx-hindsight}
Under \ref{A1} to \ref{A3}, $r_x (P_r^{\star}) \leq 1 - \frac{1}{\kappa^{\star}}$ for all $x \not \in \Xcal^\star$.
\end{lem}

Combining \Cref{thm:rx-trajectory-opt} and $\theta^{\star}_K \leq 1 - \frac{1}{\kappa^{\star}}$ from \Cref{lem:rx-hindsight}, the asymptotic linear convergence of \osgmrx~follows immediately.

\begin{coro}[Global convergence] \label{coro:rx-globalconv}
Under the same assumptions as \Cref{thm:rx-trajectory-opt}, $\theta_K^\star \leq 1 - \frac{1}{\kappa^{\star}}$ and the asymptotic complexity of {\osgmrx} to find an $\varepsilon$-optimal point is $\Ocal(\kappa^\star \log (1/\varepsilon))$, where $\kappa^\star$ is the optimal condition number defined in \eqref{eqn:optprecond}. 
\end{coro}

In fact, a slightly better convergence result can be obtained by evaluating the regret bound \eqref{eqn:rx-regret-bound} at $P = P_r^{\star}$, which we state as the following theorem. 

\begin{thm}[Refined global convergence] \label{thm:rx-globalconv}
Under \ref{A1} to \ref{A2}, \Cref{alg:osgm-rx} (\osgmrx) with \\$\eta
= \min\{ \tfrac{1}{4L^2}, \frac{\| P_r^{\star} - P_1 \|_F}{2 L \sqrt{K}} \}$ satisfies
\begin{equation} \label{eqn:ratio-global-conv}
f (x^{K + 1}) - f (x^{\star}) \leq (f (x^1) - f (x^{\star})) \big( 1 - \tfrac{1}{\kappa^{\star}} + \max\big\{\tfrac{4 L \| P^{\star}_r - P_1 \|_F}{\sqrt{K}}, \tfrac{8L^2 \| P^{\star}_r - P_1 \|_F^2}{K} \big\} \big)^K.
\end{equation}
\end{thm}

\begin{rem}
Note that \Cref{thm:rx-globalconv} has no dependence on $D$, the diameter of $\Pcal$. Therefore, {\osgmrx} can be applied even if $\Pcal = \Rbb^{n\times n}$, and there is no need to project $P$ onto $\Pcal$.
\end{rem}

The asymptotic linear convergence rate of \osgmrx~is comparable to that of preconditioned gradient descent using the universal optimal preconditioner $P_r^{\star}$. This result removes the dimension dependence from the $\Ocal(\sqrt{n}\kappa^\star \log (1/\varepsilon))$  result in \cite{kunstner2024searching}. As previously remarked, the practical convergence behavior of \osgmrx~could be even better: the linear convergence rate $\theta^{\star}_K$ is determined by the best possible choice of $P \in \mathcal{P}$ optimized for the iteration trajectory $\{x^k\}$, while the universal optimal preconditioner $P_r^{\star}$ is chosen against all possible $x$ in the domain. For convex quadratics, we have the following equivalent characterization of $P_r^\star$ through $r_x$:

\begin{prop}[Relation between $P_r^\star$ and $r_x$]\label{prop:rx-minimax}
For $f(x) = \frac{1}{2} \langle x, A x\rangle, A \in \mathbb{S}^n_{++}$, the optimal solutions to the following two problems coincide:
\begin{align}
 \min_{P \in \mathcal{P}_+} & \quad\kappa \quad
   \textrm{subject to} \quad \tfrac{1}{\kappa} I \preceq P^{1 / 2} A
   P^{1 / 2} \preceq I;  \label{eqn:minimax-sdp} \\
 \min_{P \in \mathcal{P}_+} & \max_{x \not \in \Xcal^\star } ~~ r_x(P). \label{eqn:minimax-rx} 
\end{align}
\end{prop}

Since no practical algorithm will visit every $x \in \Rbb^n$,  the trajectory-based convergence guarantee more precisely characterizes the practical performance of \osgm.

\paragraph{Unknown optimal value.}
Our method can be extended to the case where $f(x^\star)$ is unknown, but instead, a lower bound $z < f(x^\star)$ is available. In this case, we can define the auxiliary surrogate loss:
\begin{equation} \label{eqn:ratio-surrogate-loss-z}
   \rxz (P) \assign \tfrac{f (x - P \nabla f (x)) - z}{f (x) - z},
\end{equation}
obtained by replacing $f(x^\star)$ in the surrogate loss $r_x$ with lower bound $z$. Using an additional outer loop to update the lower bound $z$, the resulting algorithm (\Cref{alg:double-loop} in appendix) can achieve $\Ocal(\kappa^{\star} \log^2 (1 / \varepsilon))$ iteration complexity. 
\begin{thm}[Global convergence without knowing $f(x^\star)$, informal] \label{thm:rx-unknown-z}
Instate the same assumptions as \Cref{thm:rx-trajectory-opt} and suppose $z < f(x^\star)$ is known. There exists a variant of {\osgmrx} that finds an $\varepsilon$-optimal point in $\Ocal(\kappa^{\star} \log^2 (1 / \varepsilon))$ asymptotic complexity.
\end{thm}
The analysis for \Cref{alg:double-loop} is more involved, and we leave it to the appendix. 

\paragraph{Convex quadratics.} For strongly convex quadratics $f (x) = \frac{1}{2} \langle x, A x \rangle -
\langle b, x \rangle$, $P_r^{\star} = A^{- 1}$ gives $r_x (A^{- 1}) = 0$ for
all $x$. This implies the following superlinear convergence guarantee.

\begin{thm}[Superlinear convergence on quadratics] \label{thm:superlin-quadratic}
  For strongly convex quadratics with $\nabla^2 f(x) \equiv A \succ 0$, {\osgmrx} with $\Pcal = \Rbb^{n \times n}$ and $\eta = \frac{1}{4 L^2}$ satisfies $f (x^{K + 1}) - f (x^{\star}) \leq (f (x^1) - f (x^{\star})) (
     \tfrac{4 L^2 \| P_1 - A^{- 1} \|_F^2}{K} )^K .$
\end{thm}

\section{Gradient norm surrogate} \label{sec:gnorm-surrogate}

The second surrogate, {\textit{gradient norm surrogate}}, is  defined as follows:
\begin{equation} \label{eqn:grad-surrogate-loss}
   g_x (P) \assign \tfrac{\| \nabla f (x^+) \|}{\| \nabla f (x) \|} =
   \tfrac{\| \nabla f (x - P \nabla f (x)) \|}{\|
   \nabla f (x) \|} .
\end{equation}
Similar to the ratio surrogate, the gradient norm surrogate is defined with respect
to the contraction of the gradient norm. We assume strong convexity ($\mu > 0$) throughout this section. The gradient norm surrogate $g_x$ can be evaluated without knowing $f(x^{\star})$. However, as we will discuss later in the section, it is more challenging to establish the learnability of $g_x$, which requires the following extra assumption:
\begin{enumerate}[leftmargin=35pt,itemsep=2pt,label=\textbf{A\arabic*:},ref=\rm{\textbf{A\arabic*}},start=4]
   \item $f (x)$ has $H$-Lipschitz Hessian. $\| \nabla^2 f (x) - \nabla^2 f (y) \| \leq H \| x - y \|$ for all $x, y$.\label{A4}
\end{enumerate}

In this section, we assume a nonempty monotone oracle with respect to gradient norm $\Mcal_{\|\nabla f(x)\|, P} \not = \varnothing$.

\subsection{Surrogate loss}
Substituting $\varphi(x)$ in \Cref{thm:online-to-conv-ratio} with $\varphi (x) = \| \nabla f(x) \|$ and applying the definition of the gradient norm surrogate $g_{x}$, we obtain the following lemma.

\begin{lem}[Surrogate loss and measure] \label{lem:gx-surrogate-measure}
For all $K \geq 1$, the online scaled gradient method with nonempty monotone oracle $\Mcal_{\|\nabla f(x)\|, P} \not = \varnothing$ satisfies
  \[ \| \nabla f(x^{K+1}) \| \leq \| \nabla f(x^{1}) \|  (
     \tfrac{1}{K} \textstyle \sum_{k = 1}^K g_{x^k} (P_k)  )^K.  \]
\end{lem}

Although the gradient norm surrogate $g_x$ can be nonconvex, \Cref{prop:hx-learnability} shows that $g_x$ can be approximated by an $L$-Lipschitz continuous convex function.

\begin{prop}[Properties of $g_x$] \label{prop:gx-learnability}
  Under \ref{A1} to \ref{A4}, for any fixed $x \not \in \Xcal^\star$, the surrogate loss $g_x (P)$ defined in \eqref{eqn:grad-surrogate-loss} is $L$-Lipschitz continuous as a function in $P$ and 
  \begin{equation*}
   | g_x(P) - \hat{g}_x(P) | \leq H \| \nabla f(x) \| \|P\|^2,
  \end{equation*}
where $\hat{g}_x(P) = \Big\| \frac{\nabla f(x)}{ \|\nabla f(x)\| } - \nabla^2 f(x) P \frac{\nabla f(x)}{ \|\nabla f(x)\| } \Big\|$ is convex and $L$-Lipschitz continuous. In particular,
\[g_x(P_1) - g_x(P_2) -\langle g_x'(P_2 ), P_1 - P_2 \rangle \geq -HD^2\|\nabla f(x) \|.\]
In addition, if $x - P \nabla f(x) \not \in \Xcal^\star$, the loss $g_x(P)$ is differentiable at $P$ and its derivative takes the form
\begin{equation} \label{eqn:gx-gradient}
	\nabla g_x(P) = - \tfrac{\nabla^2 f(x- P \nabla f(x))\nabla f(x - P \nabla f(x))\nabla f(x)^\top}{\|\nabla f(x) \| \cdot \|\nabla f(x - P \nabla f(x)) \|}.
\end{equation}
\end{prop}
%\YC{If $x - P \nabla f(x) \not \in \Xcal^\star$, then $\nabla f(x^+) \neq 0$ and then $g_x$ would be differentiable. It's not clear why use subdifferential here, if we say it's a singleton, then equivalent to differentiability.}

\Cref{prop:gx-learnability} bounds the nonconvexity of $g_x$ by $\|\nabla f(x)\|$, the non-stationarity at $x$. We can still apply online learning algorithms to $g_x$ using these properties.

\subsection{Online learning algorithm}

Given Lipschitz loss functions whose nonconvexity can be bounded, online subgradient method gives the following regret guarantee.

\begin{lem}[Learnability] \label{lem:gx-learnability}
  Given \ref{A1} to \ref{A3} and the gradient norm surrogate losses $\{ g_{x^k}
  \}$, online subgradient method
\begin{equation} \label{eqn:gx-ogd}
	P_{k + 1} = \Pi_{\mathcal{P}} [P_k - \eta {g}'_{x^k} (P_k)]
\end{equation}
with stepsize $\eta = c / \sqrt{K}$ generates a sequence of scaling matrices $\{ P_k \}_{k
  \geq 2}$ such that
  \[ \textstyle \sum_{k = 1}^K g_{x^k} (P_k) - \displaystyle \min_{P \in \mathcal{P}} \textstyle \sum_{k = 1}^K g_{x^k} (P) \leq
     (\tfrac{2 D^2}{c} + \tfrac{c L^2}{2} ) \sqrt{K} + \tfrac{H
     D^2}{2} \| \nabla f (x^1) \| K. \]
  In particular, optimizing the constant $c$ suggests the stepsize $\eta = \tfrac{2 D}{L \sqrt{K}}$ and the regret bound:
\begin{equation} \label{eqn:gx-regret}
	\textstyle  \sum_{k = 1}^K g_{x^k} (P_k) - \displaystyle \min_{P \in \mathcal{P}} \textstyle \sum_{k = 1}^K
     g_{x^k} (P) \leq \rho_K \assign 2 L D  \sqrt{K} + \tfrac{H D^2}{2} \|
     \nabla f (x^1) \| K.
\end{equation}
\end{lem}

\subsection{Algorithm design and analysis}
We now state a realization of {\osgm} with the gradient norm surrogate loss $g_{x}$, denoted by {\osgmgx}.
We choose the optimality measure $\varphi$, the surrogate loss $\ell$, and the online learning algorithm $\mathcal{A}$ to be
\begin{equation*}
    \varphi(x) \assign  \|\nabla f(x)\|, \quad \ell_{x} (P) \assign g_{x} (P), \quad \mathcal{A} \assign \text{online subgradient method in \eqref{eqn:gx-ogd}},
\end{equation*}
and the monotone oracle $\Mcal$ is necessary. \Cref{alg:osgm-gx} presents the pseudocode for \osgmgx.
\begin{algorithm}[h] 
{\textbf{input} $x^1, P_1 \in \Pcal, \text{online gradient stepsize } \eta > 0, \text{nonempty } \Mcal_{\|\nabla f(x)\|,P} $}\\
\For{k =\rm{ 1, 2,...}}{
$x^{k + 1} = \Mcal_{\|\nabla f(x^k)\|, P_k} (x^{k})$\\
$P_{k+1} = \Pi_{\Pcal} [P_k - \eta {g}'_{x^k} (P_k) ]$\\
}
{\textbf{output} $x^\text{best}$ with minimum objective value}
\caption{Online scaled gradient method with gradient norm surrogate (\osgmgx) \label{alg:osgm-gx}}
\end{algorithm} 

Combining \Cref{lem:gx-surrogate-measure} and \Cref{lem:gx-learnability}, \Cref{thm:gx-trajectory-opt} characterizes the trajectory-based convergence behavior of {\osgmgx}.

\begin{thm}[Trajectory-based convergence] \label{thm:gx-trajectory-opt}
  Under \ref{A1} to \ref{A4}, \Cref{alg:osgm-gx} (\osgmgx) with \\$\eta = \frac{2D}{L \sqrt{K}}$ satisfies
  \[ \|\nabla f(x^{K+1} )\| \leq \|\nabla f(x^{1} )\| \textstyle ( \theta_K^\star + \frac{\rho_K}{K} )^K,  \]
  where $\theta^{\star}_K \assign \min_{P \in \mathcal{P}} \tfrac{1}{K}
  \sum_{k = 1}^K g_{x^k} (P)$ and $\rho_K = 2LD\sqrt{K} + \frac{HD^2}{2} \|\nabla f(x^1)\| K$ is defined in \eqref{eqn:gx-regret}.
\end{thm}

\Cref{thm:gx-trajectory-opt} itself does not necessarily yield convergence. The regret $\rho_K$ contains $\frac{HD^2}{2} \|\nabla f(x^1)\| K$ and thus is linear in $K$. One solution is to start the algorithm at a near-stationary point with sufficiently small $\|\nabla f(x^1)\|$. This strategy leads to a two-stage algorithm, and our main result is based on this strategy for brevity of exposition. 

\begin{rem}
For convex quadratics, the Lipschitz constant for Hessian is zero (i.e., $H = 0$) and the convergence of {\osgmgx} follows immediately.
\end{rem}
\subsubsection{Hindsight and global convergence}
Define $P_g^\star$ to be the gradient norm scaling matrix that solves
\[ \omega^{\star} \assign \min_{P \in \mathcal{P}} \max_x  \| I - \nabla^2 f
   (x) P \|. \]
The definition is motivated by 
\begin{align}
  \| \nabla f (x - P \nabla f (x)) \| ={} & \| \nabla f (x) - \textstyle \int_0^1
  \nabla^2 f (x - t P \nabla f (x)) P \nabla f (x) \mathd t \|
  \nonumber\\
  ={} & \| \textstyle \int_0^1 [I - \nabla^2 f (x - t P \nabla f (x)) P] \nabla f (x)
  \mathd t \| \nonumber\\
  \leq{} & [ \textstyle\int_0^1 \| I - \nabla^2 f (x - t P \nabla f (x)) P \| \mathd
  t ] \cdot \| \nabla f (x) \| \nonumber
\end{align}
and a contraction is established if $\| I - \nabla^2 f (x) P \| < 1$ for all $x$.
Define the quantity $\lambda^\star := \frac{1}{1 - \omega^\star}$. Then the following facts follow:
\begin{lem}[Hindsight] \label{lem:gx-hindsight}
  Under \ref{A1} to \ref{A3}, the followings hold:
  \begin{itemize}[leftmargin=15pt]
    \item \textit{Contraction.} $\| \nabla f (x - P^{\star}_g \nabla f (x)) \| \leq (1 - \frac{1}{\lambda^{\star}})
    \| \nabla f (x) \|$ for all $x$.
    \item \textit{Conditioning.} $\lambda^{\star} \leq  \tfrac{L}{\mu} = \kappa$.
    \item \textit{Surrogate loss bound.} $g_x(P_g^\star) \leq 1 - \frac{1}{\lambda^\star}$ for all $x \not \in \Xcal^\star$.
  \end{itemize}
\end{lem}

\begin{coro}[Global convergence] \label{coro:gx-globalconv}
 Under the same assumptions as \Cref{thm:gx-trajectory-opt}, $\theta^{\star}_K \leq 1 - \frac{1}{\lambda^\star}$ and with $\|\nabla f(x^1)\| \leq \frac{1}{H D^2\lambda^\star}$, the asymptotic complexity of {\osgmgx} to find an $\varepsilon$-critical point is $ \Ocal( 2\lambda^\star \log (1 / \varepsilon))$.

\begin{rem}
It is possible to sharpen \Cref{coro:gx-globalconv} using a more fine-grained analysis: the nonconvexity will vanish as the algorithm converges, and $\Ocal(\lambda^\star \log (1 / \varepsilon))$ complexity still holds asymptotically.
\end{rem}
\begin{rem}
{\osgmgx} can also output an $\varepsilon$-optimal solution due to the relation $f(x) - f(x^\star) \leq \frac{1}{2\mu}\|\nabla f(x)\|^2$ from strong convexity.
\end{rem}
	
\end{coro}
\section{Hypergradient surrogate} \label{sec:hypergrad-surrogate}

The last surrogate loss, {\textit{hypergradient surrogate}}, is defined as follows:
\begin{equation} \label{eqn:hx-definition}
	h_x (P) \assign \tfrac{f (x^+) - f (x)}{\| \nabla f (x) \|^2} = \tfrac{f (x - P \nabla f(x)) - f (x)}{\| \nabla f (x) \|^2} .
\end{equation}
The name hypergradient comes from {\cite{gunes2018online}}, by which the hypergradient descent heuristic improves the convergence of
gradient-based methods. Unlike the ratio surrogate $r_x$ or the gradient norm surrogate $g_x$, the hypergradient surrogate $h_x$ itself is not directly derived from $\varphi (x)$ using a telescopic product, but instead motivated by the descent lemma:
\[ f \big( x - \tfrac{1}{L} \nabla f (x) \big) - f (x) \leq - \tfrac{1}{2
   L} \| \nabla f (x) \|^2. \]
Dividing both sides of the inequality by $\| \nabla f (x) \|^2$ for $x \not \in \Xcal^\star$ gives $h_x$. The descent lemma does not depend on the strong convexity coefficient $\mu$, so the hypergradient surrogate $h_x$ applies to general convex (non-strongly convex) optimization problems. Throughout this section, we assume a nonempty monotone oracle with respect to function value gap $\Mcal_{f(x) - f(x^\star), P} \not = \varnothing$.

\subsection{Surrogate loss}

To analyze the hypergradient surrogate $h_x$, we must connect it with a measure of convergence. \Cref{lem:hx-surrogate-measure} presents this relation.
\begin{lem}[Surrogate loss and measure] \label{lem:hx-surrogate-measure}
  Under \ref{A1}, \ref{A2}, for all $K \geq 1$, the online scaled gradient method with nonempty monotone oracle $\mathcal{M}_{f(x) - f(x^\star), P} \not = \varnothing$ satisfies:
  \begin{itemize}[leftmargin=15pt,itemsep=-10pt, topsep=5pt]
    \item If $\mu > 0$, then
\begin{equation}
	f (x^{K + 1}) - f (x^{\star}) \leq (f (x^1) - f (x^{\star}))
       ( 1 - 2 \mu \max \{ \tfrac{1}{K} \textstyle \sum_{k = 1}^K -h_{x^k} (P_k),
       0 \} )^K \label{eqn:hx-strong-cvx}.
\end{equation}
    \item If $\mu \geq  0$, then 
	\begin{align}
		\min_{1 \leq k \leq K}  \| \nabla f (x^k) \|^2 \leq{} & \tfrac{f (x^1)
       - f (x^{\star})}{K} \cdot\tfrac{1}{\max \{ \frac{1}{K} \sum_{k = 1}^K -h_{x^k}
       (P_k), 0 \}} \label{eqn:hx-cvx-gnorm}, \\
       f (x^{K + 1}) - f (x^{\star}) \leq{} & \tfrac{\Delta^2}{K} \cdot \tfrac{1}{\max
       \{  \frac{1}{K} \sum_{k = 1}^K -h_{x^k} (P_k), 0 \}}, \label{eqn:hx-cvx-fval}
	\end{align}
  \end{itemize}
  ~~\quad where $\Delta = \max_{x \in \Lcal_{f(x^1)}} \min_{x^{\star} \in
  \mathcal{X}^{\star}} \| x - x^{\star} \|$.
\end{lem}

\begin{rem}
The $\max \{\cdot, 0\}$ terms arise from the monotone oracle.
Here, we slightly abuse the notation: if the denominator in \eqref{eqn:hx-cvx-gnorm} or \eqref{eqn:hx-cvx-fval} is 0, the bound simplifies to a trivial bound, with the right-hand side being infinity. 
Note that the surrogate loss appears in the denominator of the sublinear convergence rate, which differs from the previous analyses.
\end{rem}

Now we establish the properties of the hypergradient surrogate $h_x$.
\begin{prop}[Properties of $h_x$] \label{prop:hx-learnability}
  Under \ref{A1} to \ref{A3}, for any fixed $x \not \in \Xcal^\star $, the surrogate loss $h_x (P)$ defined in \eqref{eqn:hx-definition} is convex
  and $(LD + 1)$-Lipschitz continuous as a function in $P$. In addition, the derivative of $h_x$ takes the form
\begin{equation} \label{eqn:hx-gradient}
\nabla h_{x} (P) = - \tfrac{\nabla f (x - P \nabla f(x)) \nabla f (x)^{\top}}{\|\nabla f(x)\|^2}.
\end{equation}
\end{prop}

\subsection{Online learning algorithm}

Given convex and Lipschitz-continuous losses, online gradient descent gives the
following regret guarantee.

\begin{lem}[Learnability] \label{lem:hx-learnability}
  Given \ref{A1} to \ref{A3} and the hypergradient surrogate losses $\{ h_{x^k} \}$, online gradient descent
\begin{equation} \label{eqn:hx-ogd}
	P_{k + 1} = \Pi_{\mathcal{P}} [P_k - \eta \nabla h_{x^k} (P_k)]
\end{equation}
  with stepsize $\eta = c / \sqrt{K}$ generates a sequence of scaling matrices $\{ P_k \}_{k
  \geq 2}$ such that
  \[ \textstyle \sum_{k = 1}^K h_{x^k} (P_k) - \displaystyle \min_{P \in \mathcal{P}} \textstyle \sum_{k = 1}^K h_{x^k}
     (P) \leq (\tfrac{2D^2}{c} + \tfrac{c (LD + 1)^2}{2}) \sqrt{K}. \]
     In particular, optimizing the constant $c$ suggests the stepsize $\eta = \frac{2D}{(LD + 1) \sqrt{K}}$ and the regret bound:
\begin{equation} \label{eqn:hx-regret}
	\textstyle \sum_{k = 1}^K h_{x^k} (P_k) - \displaystyle \min_{P \in \Pcal} \textstyle\sum_{k = 1}^K h_{x^k}
     (P) \leq \rho_K := 2D(LD + 1) \sqrt{K}.
\end{equation}
\end{lem}

\subsection{Algorithm design and analysis}
We now state a realization of {\osgm} with the hypergradient surrogate loss $h_{x}$, denoted by {\osgmhx}. 
We choose the optimality measure $\varphi$, the surrogate loss $\ell$, and the online learning algorithm $\mathcal{A}$ to be
\begin{equation*}
    \varphi(x) \assign f(x) - f(x^\star) \text{ or } \|\nabla f(x)\|, \quad \ell_{x} (P) \assign h_{x} (P), \quad \mathcal{A} \assign \text{online gradient descent in \eqref{eqn:hx-ogd}},
\end{equation*}
and the monotone oracle $\Mcal$ is necessary. \Cref{alg:osgm-hx} presents the pseudocode for {\osgmhx}.
\begin{algorithm}[h]
{\textbf{input} $x^1, P_1 \in \Pcal, \text{online gradient stepsize } \eta > 0, \text{nonempty oracle }\Mcal_{f(x) - f(x^\star),P} $}\\
\For{k =\rm{ 1, 2,...}}{
$x^{k + 1} = \Mcal_{f(x) - f(x^\star), P_k}(x^{k})$\\
$P_{k+1} = \Pi_{\Pcal} [P_k - \eta \nabla h_{x^k} (P_k) ]$\\
}
{\textbf{output} $x^\text{best}$ with minimum objective value}
\caption{Online scaled gradient method with hypergradient surrogate (\osgmhx) \label{alg:osgm-hx}}
\end{algorithm} 

Combining \Cref{lem:hx-surrogate-measure} and \Cref{lem:hx-learnability}, \Cref{thm:hx-trajectory-opt} characterizes the trajectory-based convergence behavior of {\osgmhx}.
\begin{thm}[Trajectory-based convergence] \label{thm:hx-trajectory-opt}
  Under \ref{A1} to \ref{A3}, \Cref{alg:osgm-hx} (\osgmhx) with \\$\eta = \frac{2D}{(LD + 1)\sqrt{K}}$ satisfies
  \begin{itemize}[leftmargin=15pt,itemsep=-10pt, topsep=5pt]
    \item If $\mu > 0$, then
    \[ f (x^{K + 1}) - f (x^{\star}) \leq (f (x^1) - f (x^{\star}))
       ( 1 - 2 \mu \max \{-\theta^{\star}_K - \tfrac{\rho_K}{K}, 0
       \} )^K. \]
    \item If $\mu \geq  0$, then 
	\begin{align*}
		\min_{1 \leq k \leq K}  \| \nabla f (x^k) \|^2 \leq{} & \tfrac{f (x^1)
       - f (x^{\star})}{K} \cdot \tfrac{1}{\max \{ -\theta^{\star}_K - \frac{\rho_K}{K}, 0 \}},\\
       f (x^{K + 1}) - f (x^{\star}) \leq{} & \tfrac{\Delta^2}{K} \cdot \tfrac{1}{\max \{ -\theta^{\star}_K - \frac{\rho_K}{K}, 0 \}},
	\end{align*}
  \end{itemize}
  where $\theta^{\star}_K \assign \min_{P \in \mathcal{P}} \tfrac{1}{K}
  \sum_{k = 1}^K h_{x^k} (P)$, $\Delta$ is defined in \Cref{lem:hx-surrogate-measure} and $\rho_K = 2D(LD + 1) \sqrt{K}$ is defined in \eqref{eqn:hx-regret}.
\end{thm}

\subsubsection{Hindsight and global convergence}

Define $P^{\star}_h$ to be the hypergradient scaling matrix that
solves
\[ \gamma^{\star} \assign \max_{P \in \mathcal{P}} \min_{~x \in \Lcal_{f(x^1)} \backslash
   \mathcal{X}^{\star}} - h_x(P) = \tfrac{f(x) - f(x - P\nabla f(x)}{\|\nabla f(x)\|^2}. \]
    
Intuitively, $\gamma^\star$ maximizes the function value progress with respect to the gradient norm and can be interpreted as the inverse of the local Lipschitz smoothness constant. The descent lemma gives a lower bound on $\gamma^\star$.

\begin{lem}[Hindsight] \label{lem:hx-hindsight}
  Under \ref{A1} to \ref{A3}, $-h_x(\Pcal_h^\star) \geq \gamma^\star \geq  \tfrac{1}{2 L}$ for all $x \not \in \Xcal^\star$.
\end{lem}
\begin{coro}[Global convergence]  
\label{coro:hx-globalconv}
Under the same assumptions as \Cref{thm:hx-trajectory-opt}, 
$\theta_K^\star \leq -\gamma^\star$, and
\begin{itemize}[leftmargin=15pt]
	\item If $\mu > 0$, then the asymptotic complexity of {\osgmhx} to find an $\varepsilon$-optimal point is $\Ocal(\frac{1}{2\mu\gamma^\star} \log (1/\varepsilon))$.
	\item If $\mu = 0$, then the asymptotic complexity of {\osgmhx} to find an $\varepsilon$-optimal point is $\Ocal(\frac{1}{\gamma^\star\varepsilon})$.
\end{itemize}
\end{coro}

\begin{rem}
Given $-\theta_K^\star \geq  \tfrac{1}{2 L}$, the complexity of {\osgmhx} is no worse than vanilla gradient descent and can provide acceleration if $-\theta_K^\star > \tfrac{1}{2 L}$. Our results for the first time show that the hypergradient heuristic, when combined with a monotone oracle $\Mcal$,	provably accelerates gradient descent.
\end{rem}

\section{Practical considerations} \label{sec:practical}

This section considers practical aspects of implementing {\osgm}. 

\subsection{Gradient evaluations for $r_x$ and $h_x$ using simple comparison oracle}
A first look at \Cref{alg:osgm-rx} and \Cref{alg:osgm-hx} suggests an additional gradient evaluation at every iteration. Both $\nabla r_x$ and $\nabla h_x$ need to evaluate two gradients in every iteration.  However, with $x^{k + 1/2} := x^k - P_k \nabla f(x^k)$ and simple comparison (\Cref{sec:monotone-oracle}) as the monotone oracle, the gradient $\nabla r_{x^k} (P)$ in \Cref{alg:osgm-rx} can be expressed as
\begin{equation} \label{eqn:grad-ratio-xk}
   \nabla r_{x^k} (P) = - \tfrac{\nabla f (x^{k+1/2}) \nabla f (x^k)^{\top}}{f (x^k) - f (x^{\star})}.
\end{equation} 
If $x^{k+1} = x^{k + 1/2}$, then $\nabla f (x^{k+1/2})$  can be reused in the next iteration; if $x^{k+1} = x^{k}$, then $\nabla f (x^{k})$  can be reused. Therefore, simple comparison oracle ensures that the total number of gradient evaluations in \osgmrx~is the same as that of gradient descent but simply requires an additional cache to store $\nabla f(x^{k+1})$. Regarding $\nabla g_x$, its implementation needs a Hessian-gradient product, which can be efficiently computed in practice.

\subsection{Efficient scaling matrix updates}
The subset $\mathcal{P}$ can be chosen to have a simple structure, such as diagonal matrices or sparse matrices with some prespecified sparsity pattern. Then the scaling matrix $P_{k+1}$ can be efficiently updated in the cost of $\mathcal O (\texttt{supp}(\mathcal P))$ since it suffices to compute the non-zero entries of the sparsity pattern in $\mathcal{P}$.
For example, if $\mathcal{P}$ is the set of diagonal matrices, then \eqref{eqn:grad-ratio-xk} simplifies to $\nabla r_{x^k} (P) = - \tfrac{\nabla f (x^{k} - P_k \nabla f(x^k)) \odot \nabla f (x^k)}{f (x^k) - f (x^{\star})}$, where $\odot$ denotes the element-wise product. 
A simpler structure in $\mathcal{P}$ provides more efficient scaling matrix updates. However, more freedom in $\mathcal{P}$ may provide smaller $\theta^{\star}_K$, enhancing the convergence of \osgm.

\subsection{Choice of candidate set of scaling matrices $\mathcal{P}$}

We propose two heuristics for choosing a subset $\mathcal{P}$ of
sparse matrices. Sparsity refers to either entries sparsity or spectrum sparsity. We assume that some Hessian matrix $\nabla^2 f (x) = A \succ 0$ is
known.

\begin{itemize}[leftmargin=10pt]
	\item \textit{Nonzero sparsity pattern.}\\
	A preconditioner can be viewed as a cutting plane in the
difference of extremal spectrum {\cite{gao2023scalable}}. Let $v_{\min}$ and
$v_{\max}$ be two extremal eigenvectors of $A$. Then $| v_{\max}
v_{\max}^{\top} - v_{\min} v_{\min}^{\top} |$, an $n \times n$ grid with
nonnegative entries, serves as a score function for the most critical sparsity
pattern. The large-magnitude entries in $| v_{\max}
v_{\max}^{\top} - v_{\min} v_{\min}^{\top} |$ strongly affect the conditioning of the system.
\item \textit{Spectral sparsity (low rank).}\\
It is common to consider diagonal plus low-rank preconditioners, and
randomized preconditioners have proved to be very efficient
{\cite{frangella2023randomized}}. Given a low-rank matrix $M$, we can
parameterize $\mathcal{P}= \{ \tmop{diag} (d) + \alpha M : (d, \alpha) \in
\mathbb{R}^{n + 1} \}$ to be the linear combination between diagonal
matrices and $M$.

\end{itemize}

\subsection{Choice of online learning algorithm $\Acal$}
Our convergence analyses show that a good online learning algorithm $\Acal$ benefits convergence of {\osgm}. For simplicity, the simplest possible online learning algorithms are adopted in the theoretical analysis. However, {\osgm} is compatible with more advanced online learning algorithms such as {\adagrad}. Using advanced online algorithms often greatly improves the robustness and practical performance of {\osgm}. In particular, our results on the hypergradient surrogate loss provide new insights into improving the hypergradient descent heuristics.

\section{Numerical experiments} \label{sec:exp}

In this section, we conduct experiments to show the performance of
online scaled gradient methods. We test on standard strongly convex
optimization problems, including least squares and regularized logistic
regression.

\subsection{Experiment setup}

\paragraph{Synthetic data.}For the least squares problem $f(x) = \tfrac{1}{2} \| A x - b \|^2$,  $A \in \mathbb{R}^{n
\times n} = C D C^\top + \sigma  I$ with $C$ is element-wise generated from $0.01\times \mathcal{N} (0, 1)$ and an identity matrix $I$ is added to it; $D$ is a diagonal matrix with $\mathcal{U} [0, 1]^n$ diagonals; $b \in \mathbb{R}^m$ is generated from $\mathcal{U} [0, 1]^n$. 

\paragraph{Real data.}We use datasets from \texttt{LIBSVM}
{\cite{chang2011libsvm}} for support vector machine (SVM) problem $f (x) = \frac{1}{m}
\sum_{i = 1}^m f_i (x) + \frac{\lambda}{2} \| x \|^2$, where $f_i$ is the squared hinge loss \cite{lee2001ssvm}. We set $\lambda = 5 / n$.

\paragraph{Benchmark algorithms.} Eight algorithms are compared:
\begin{itemize}[leftmargin=10pt]
  \item ({\gd}) Gradient descent with $1 / L$ stepsize.
  \item ({\optdgd}) Gradient descent with the universal optimal diagonal preconditioner \cite{qu2024optimal, gao2023scalable}.
  \item ({\osgmrx}) Online scaled gradient method with ratio surrogate.
  \item ({\osgmgx}) Online scaled gradient method with gradient norm surrogate.
  \item ({\osgmhx}) Online scaled gradient method with hypergradient surrogate.
  \item ({\adagrad}) Adaptive (sub)gradient method \cite{duchi2011adaptive}.
  \item ({\agd}) Accelerated gradient descent for general convex problems \cite{d2021acceleration}. 
  \item ({\sagd}) Accelerated gradient descent for strongly convex problems \cite{d2021acceleration}.
\end{itemize}

{\optdgd} and {\osgmgx} are only tested on problems with fixed Hessian.

\paragraph{Algorithm configurations.} We configure the algorithms as follows.
\begin{enumerate}[leftmargin=15pt, label=\textbf{\arabic*)}]
\item \textit{Dataset generation}. For synthetic data, we pick $n = 100$ and $\sigma \in \{10^{-4}, 10^{-3}, 10^{-2}, 10^{-1}\}$.
  \item \textit{Initial point}. Initial points $x^1$ for all the algorithms are generated from standard normal $\Ncal(0, I_n)$ and scaled to have unit length. Initial scaling matrix $P_0 = 0$.
  
  \item \textit{Maximum iteration}. The maximum iteration is set to $K = 10000$.
  \item \textit{Stopping criterion}. Algorithm stops when $\|\nabla f(x^k)\| \leq 10^{-10}$.
 
  \item \textit{Stepsize configuration}. ({\adagrad}) uses the optimal stepsize
  among $\{ 10^{- 3}, 10^{- 2}, 10^{- 1}, 1, 10 \}$.
  
  \item \textit{Monotone oracle}. All {\osgm} methods use simple comparison (\Cref{sec:monotone-oracle}) as the monotone oracle. 
  
  \item \textit{Online learning algorithm}. All {\osgm} methods use $\Acal ={}$({\adagrad}) with the optimal stepsize among $\{ 10^{- 3}, 10^{- 2}, 10^{- 1}, 1, 10 \}$. 
  
  \item \textit{Choice of candidate scaling matrix $\mathcal{P}$}. We choose $\mathcal{P} = \Rbb^{n\times n}$ in \Cref{exp:superlin} and $\mathcal{P}$ as the
  set of diagonal matrices in the rest of the experiments.
  
  \item \textit{Knowledge of optimal value}. When the exact optimal value is unknown, {\osgmrx} uses the auxiliary surrogate $\rxz$ in \eqref{eqn:ratio-surrogate-loss-z}. But we allow $z$ to be an arbitrary guess of $f(x^\star)$, and if $z > f(x^k)$, we heuristically adjust $z \leftarrow f(x^k) - \min \{5(z - f(x^k)), 1\}$ to update the lower bound. This strategy is not theoretically justified but performs well in practice.
\end{enumerate}

\subsection{Toy example: near diagonal convex quadratic}
\begin{figure}[h]
\centering
\includegraphics[scale=0.19]{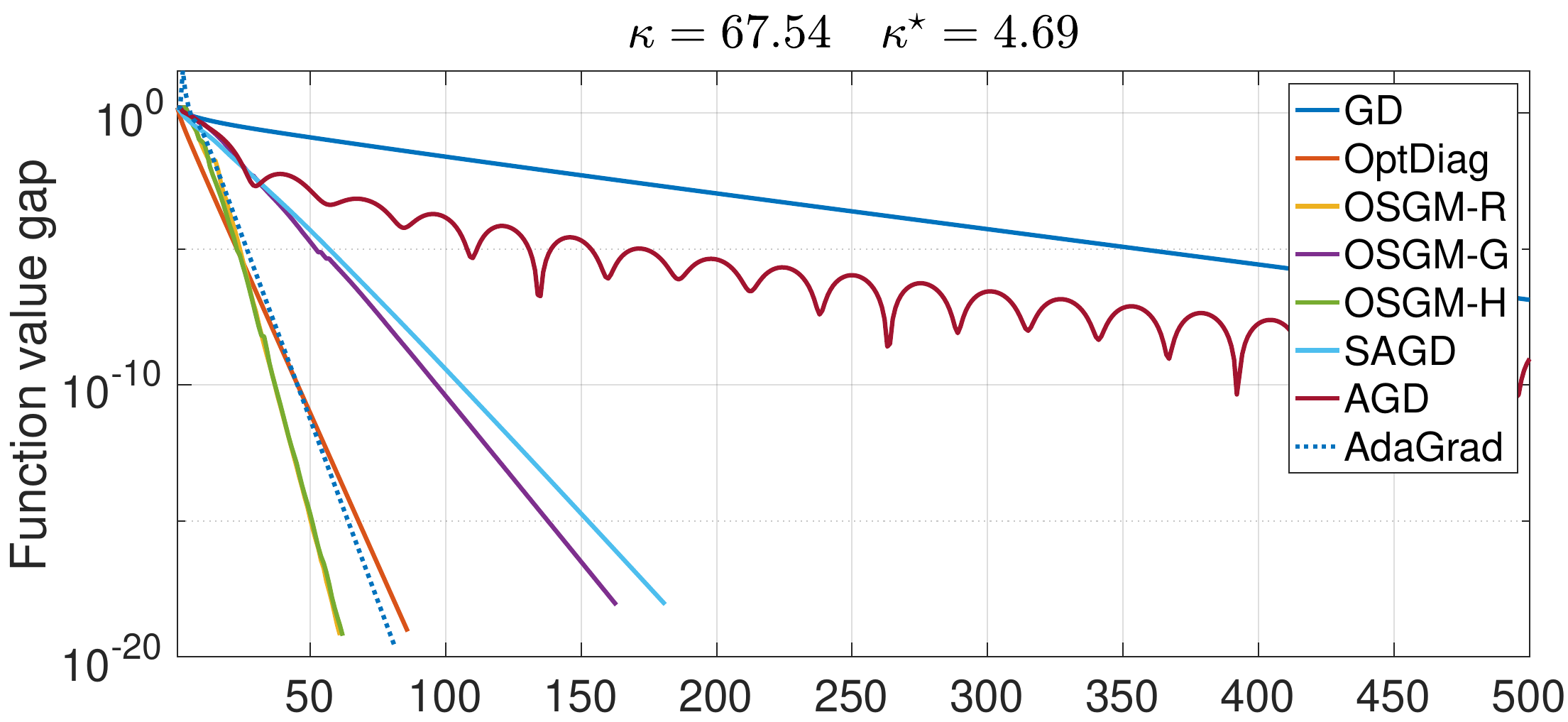}\qquad\includegraphics[scale=0.21]{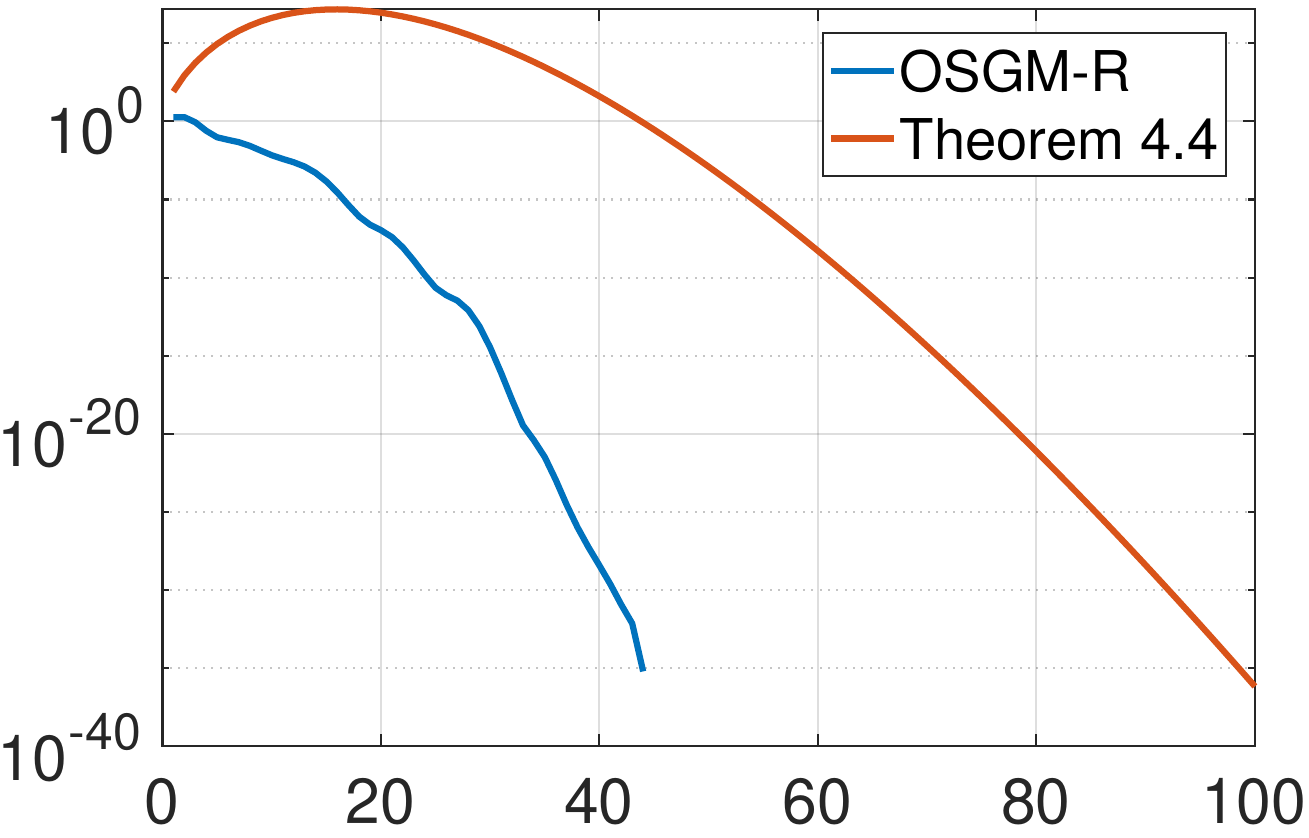}
\caption{Left: comparison of benchmark algorithms on toy quadratic problem. \label{fig:toy} Right: superlinear convergence of {\osgmrx} on convex quadratics. x-axis: iteration count.}
\end{figure}

This section verifies the convergence behavior of {\osgm} on a toy least squares problem with near diagonal Hessian. The problem has $\kappa \approx 68$ and $\kappa^\star \approx 4.7 < \sqrt{\kappa}$. Theory predicts that {\osgm} should outperform {\sagd} asymptotically.  \Cref{fig:toy} (left) illustrates the performance of the eight tested algorithms, with {\osgmrx} and {\osgmhx} showing the most competitive performance. In particular, the linear convergence rates (slope) of three {\osgm} algorithms are better than that of {\sagd}, which aligns with our theory. Moreover, {\osgmrx} and {\osgmhx} both converge faster than gradient descent using the universal diagonal preconditioner ({\optdgd}). This is also consistent with our theory and suggests that we can still gain from being adaptive, even on a convex quadratic with fixed curvature. Notably, {\adagrad} also achieves competitive performance on this problem.

\subsection{Superlinear convergence on quadratics} \label{exp:superlin}
This section verifies the \Cref{thm:superlin-quadratic}, the superlinear convergence behavior of {\osgmrx}. We take $\mathcal{P} = \Rbb^{n\times n}$ and plot the theoretical bound in \Cref{thm:superlin-quadratic} as well as the true performance of {\osgmrx}. \Cref{fig:toy} (right) indicates that {\osgmrx} exhibits superlinear convergence and validates our theory.

\subsection{More comparison between the algorithms}

\begin{figure}[h]
\centering
\includegraphics[scale=0.21]{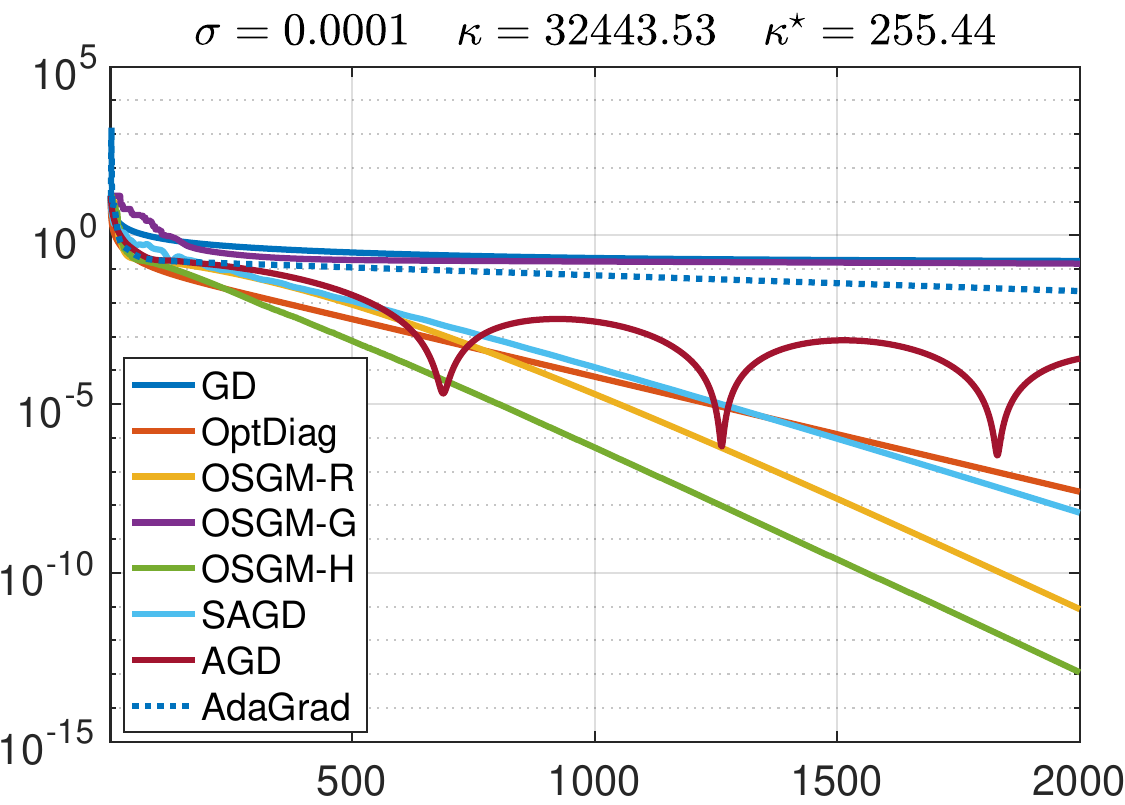}
\includegraphics[scale=0.21]{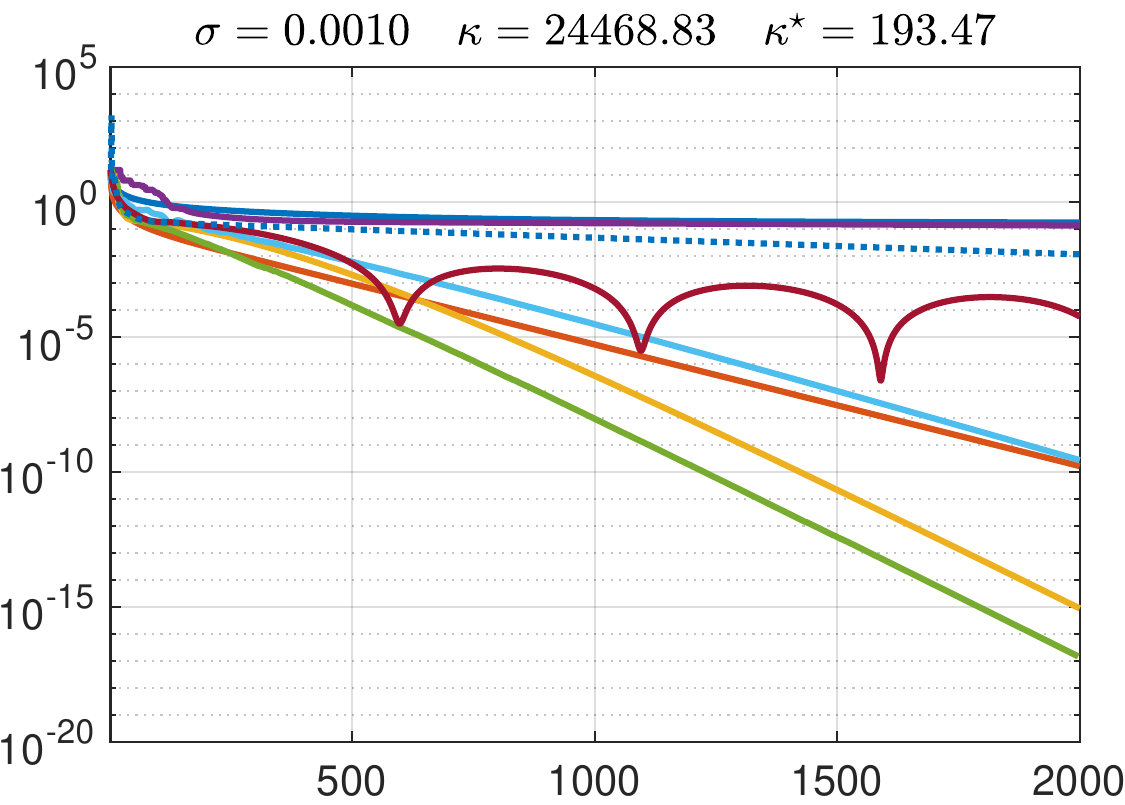}
\includegraphics[scale=0.21]{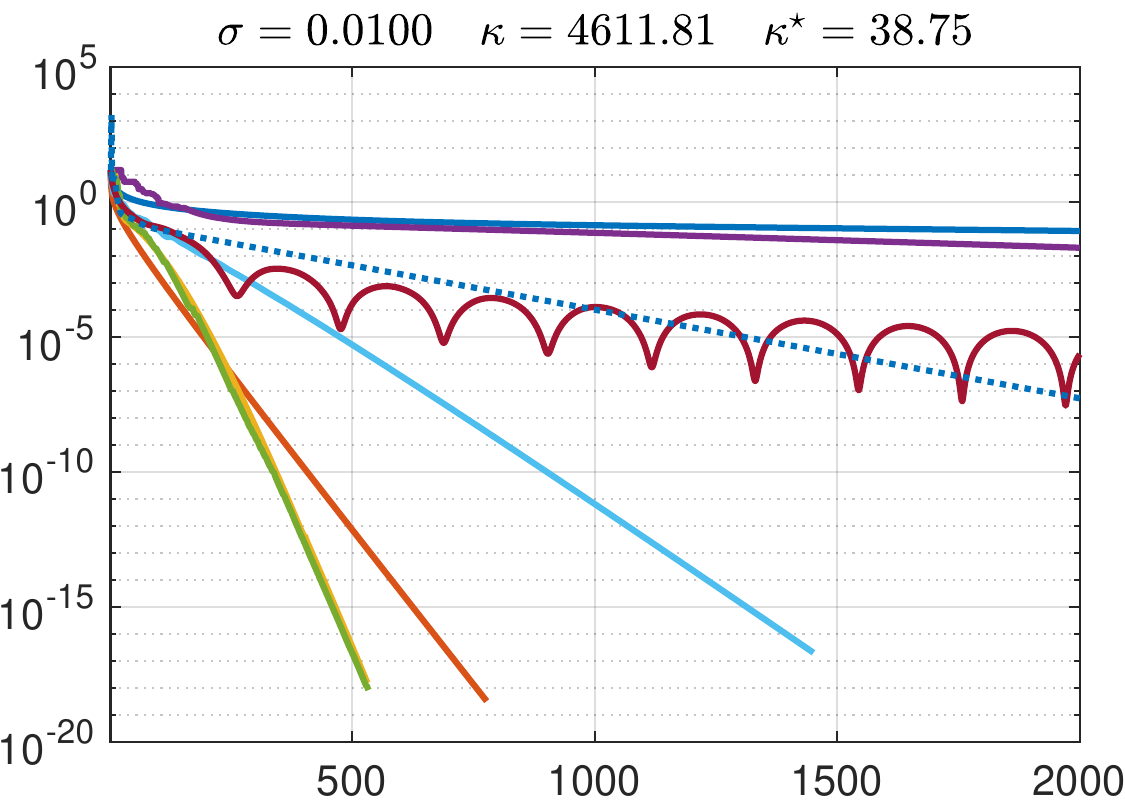}
\includegraphics[scale=0.21]{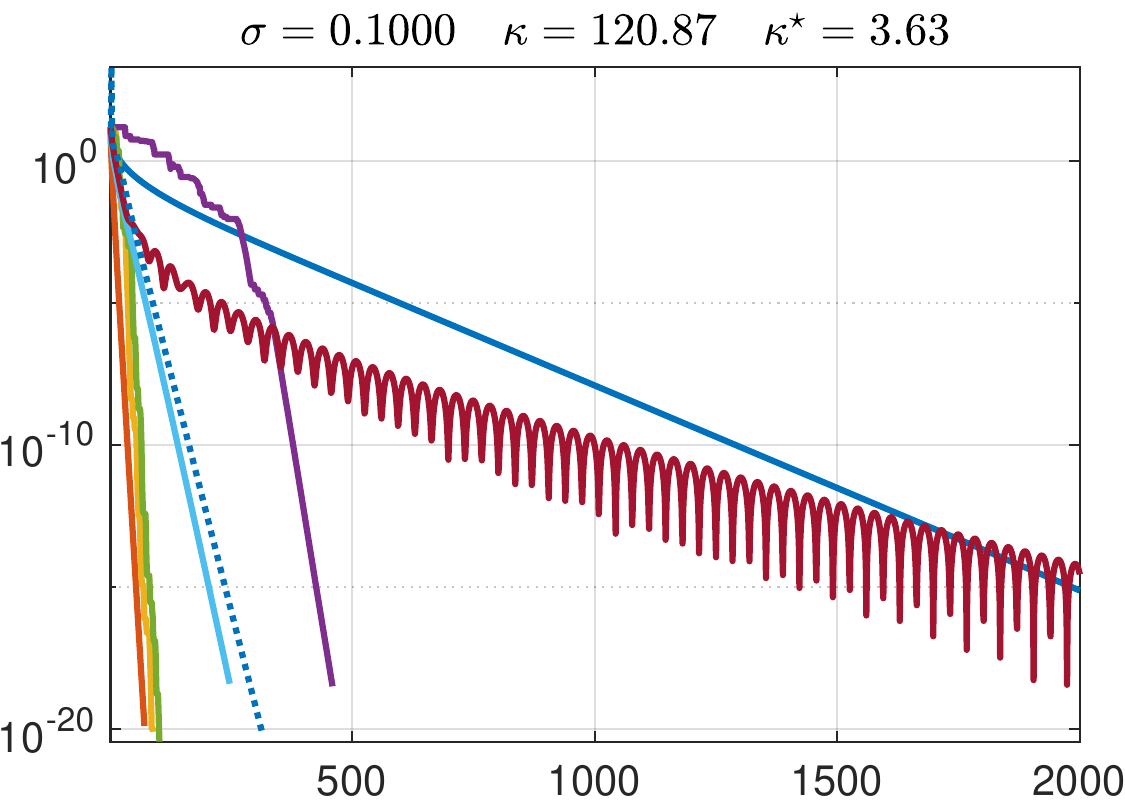}
\caption{Function value gap on least squares problem with $\sigma \in \{
  10^{- 4}, 10^{- 3}, 10^{- 2}, 10^{- 1} \}$ \label{fig:ls}}
\end{figure}

This section compares different algorithms on the aforementioned datasets. \Cref{fig:ls} shows the results on synthetic least squares problems with different condition numbers.

\Cref{fig:ls} suggests that when $\kappa^\star \ll \sqrt{\kappa}$, {\osgm} tends to outperform accelerated methods. On the other hand, if $\kappa^\star > \sqrt{\kappa}$, {\osgm} is often less competitive compared to {\sagd} on quadratics. \\

\begin{figure}[h]
\centering
\includegraphics[scale=0.21]{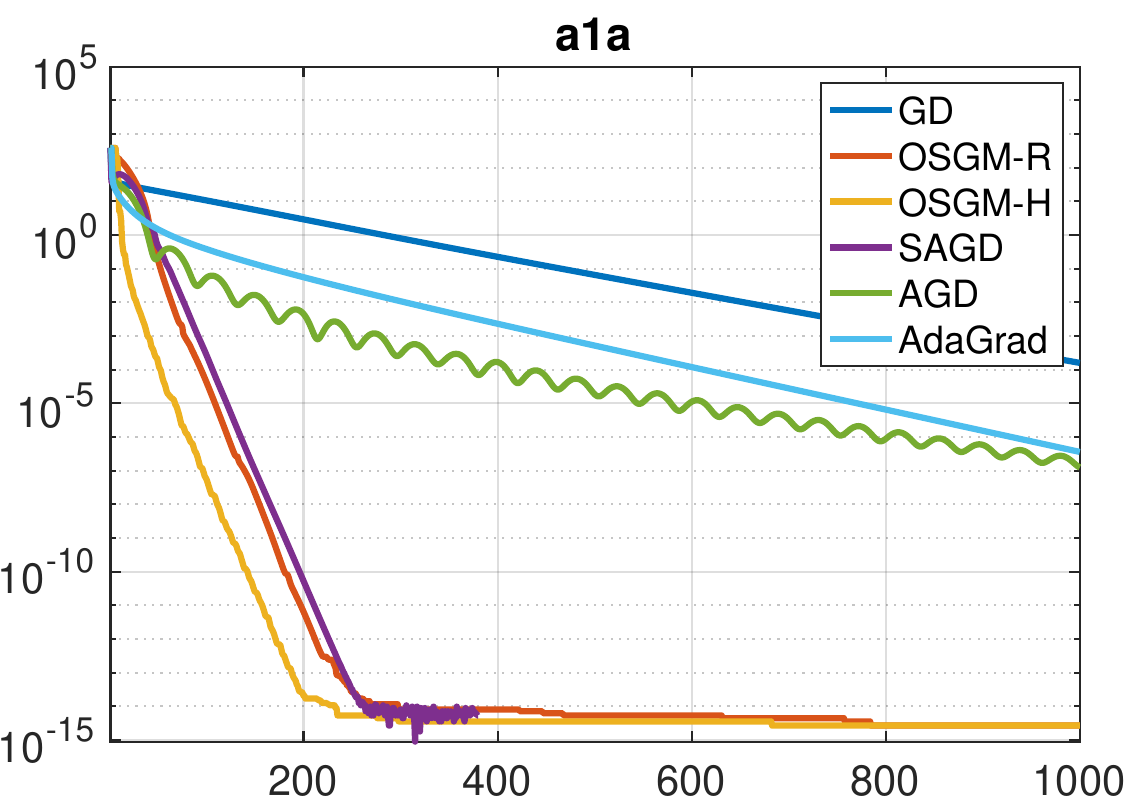}
\includegraphics[scale=0.21]{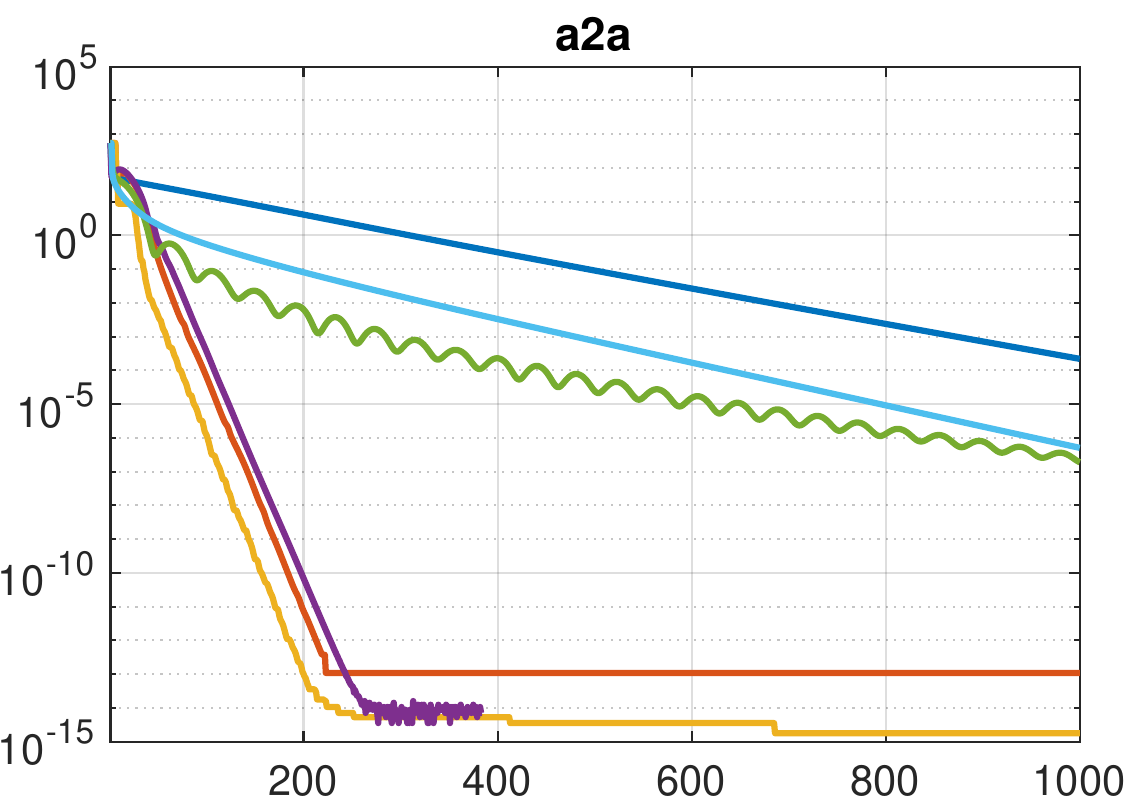}
\includegraphics[scale=0.21]{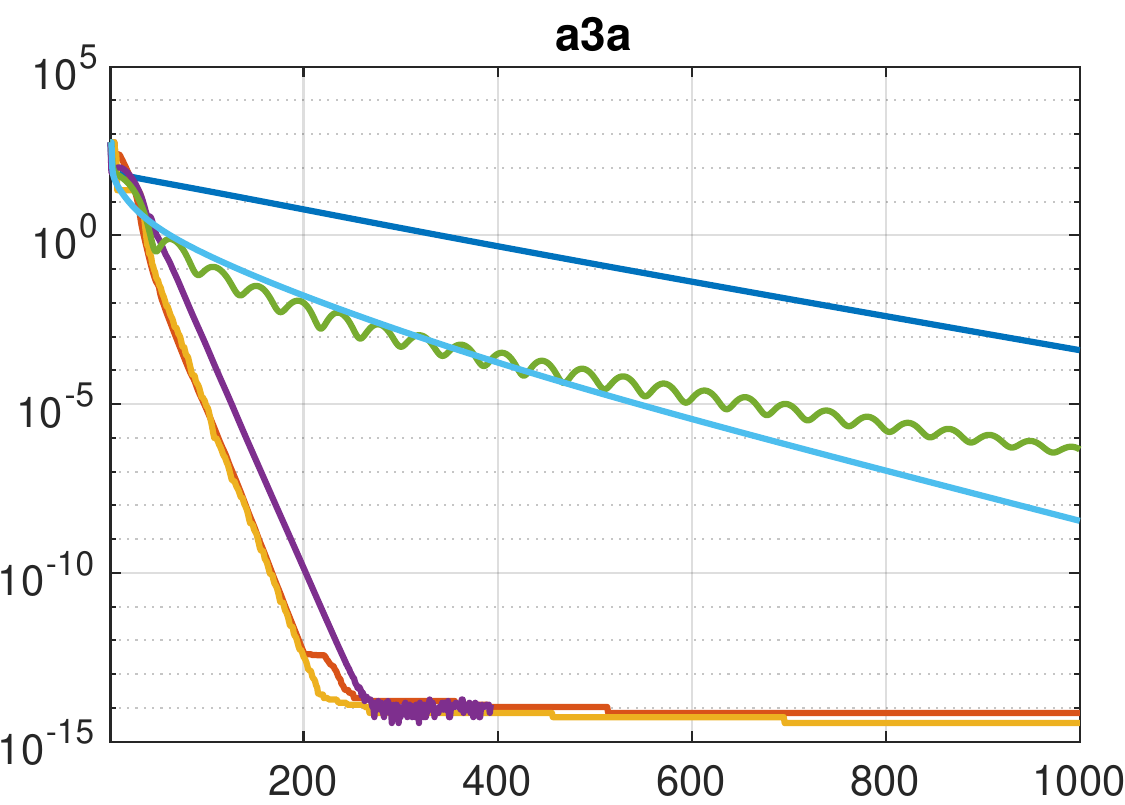}
\includegraphics[scale=0.21]{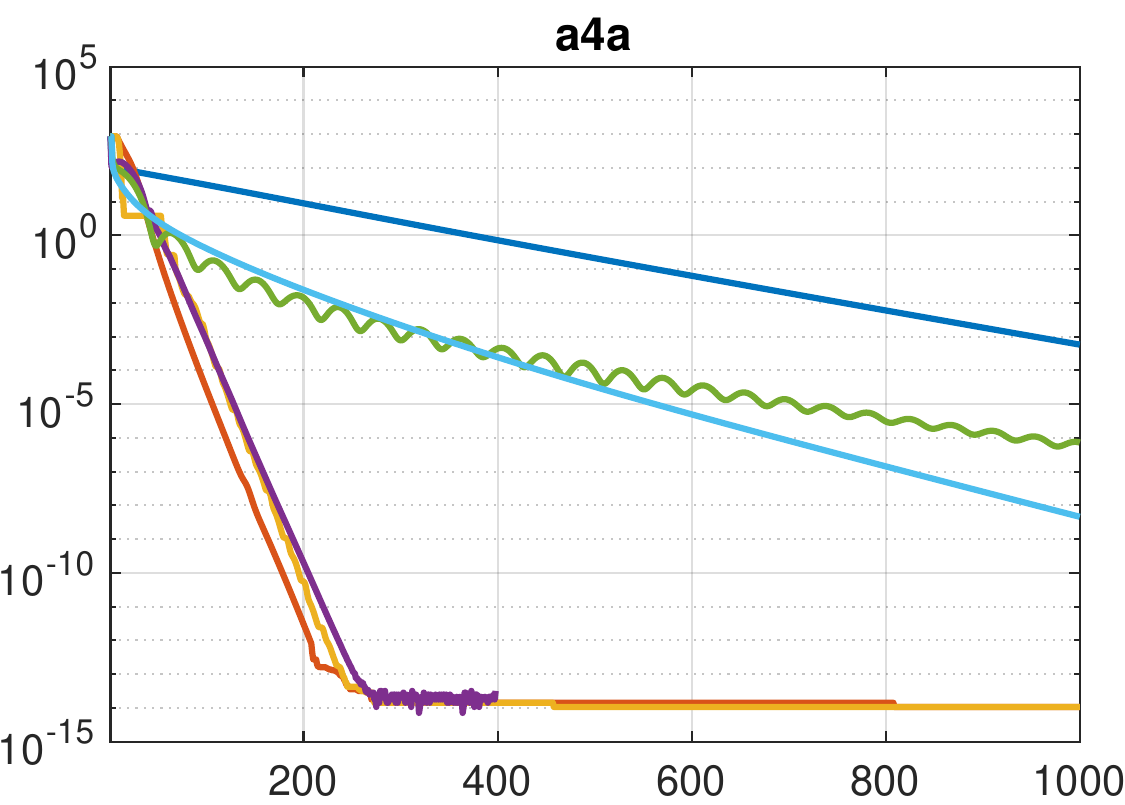}
\caption{Function value gap on the support vector machine problems \label{fig:svm}}
\end{figure}

\textbf{Figure \ref{fig:svm}} shows the results on the SVM problems from \texttt{LIBSVM}, and we observe similar competitive performance of {\osgmrx} and {\osgmhx}.

\section{Conclusions and future directions}

In this paper, we discuss \osgm, a general framework that allows online convex optimization algorithms to provably accelerate gradient-based algorithms. Our framework achieves a strong trajectory-based convergence guarantee and explains the success of the popular hypergradient descent heuristic. Future directions include extending the results to accelerated gradient descent, stochastic gradient descent, nonconvex, nonsmooth, and constrained optimization, and to other iterative algorithms where a scaled update affects the algorithm performance. 
\section*{Acknowledgement}
We appreciate the constructive feedback from Professor Qi Deng from Shanghai Jiaotong University and Zikai Xiong from MIT.

\renewcommand \thepart{}
\renewcommand \partname{}

\bibliography{ref.bib}
\bibliographystyle{plain}

\doparttoc
\faketableofcontents
\part{}

\newpage
\appendix

\addcontentsline{toc}{section}{Appendix}
\part{Appendix} 
\parttoc

\newpage
\section{Proof of results in Section \ref{sec:online}}

\subsection{Proof of Theorem \ref{thm:online-to-conv-ratio}}

Since the measure $\varphi$ is non-negative. Applying arithmetic-geometric mean inequality 
\[ \textstyle \big(\prod_{k=1}^K a_k\big)^{1/K} \leq \frac{1}{K}\sum_{k=1}^K a_k\]
completes the proof.

\section{Proof of results in Section \ref{sec:ratio-surrogate}}

\subsection{Auxiliary results}

\begin{lem}[\cite{orabona2019modern}] \label{lem:auxi-smooth}
  Let $r (P)$ be a $\tau$-smooth function with $\min_{P \in \mathcal{P}} r (P) \geq 0$. Then $r (P) \geq \tfrac{1}{2 \tau} \| \nabla r (P) \|^2$ for all $P \in \mathcal{P}$.
\end{lem}

\begin{lem} \label{lem:smooth-online-learnability}
  Given a family of non-negative, convex, and $\tau$-smooth losses $\{\rk\}$, online gradient descent
  \begin{equation}\label{eqn:online-pgd}
    P_{k + 1} = \Pi_{\mathcal{P}} [P_k - \eta \nabla \rk (P_k)]
  \end{equation}
  with stepsize $\eta \leq 1/(2 \tau)$ generates a sequence of scaling matrices $\{P_k\}_{k \geq 2}$ such that 
  \begin{equation} \label{eqn:r-regret-bound-z}
      \textstyle \sum_{k = 1}^K \rk (P_k) - \textstyle \sum_{k = 1}^K \rk (P) \leq \tfrac{1}{\eta} \| P - P_1 \|^2_F + 2 \tau \eta \textstyle \sum_{k = 1}^K \rk (P) \quad \text{for any } P \in \mathcal{P}.
  \end{equation}
\end{lem}

\begin{proof}
  The proof follows the standard proof of the $L^{\star}$ regret bound \cite{orabona2019modern} in online convex optimization but is tailored to our settings. 
For any $P \in \mathcal{P}$, we have
\begin{align}
  \| P_{k + 1} - P \|^2_F & ={} \| \Pi_{\mathcal{P}} [P_k - \eta \nabla \rk
  (P_k)] - P\|^2_F \nonumber \\
  & \leq{} \| P_k - P - \eta \nabla \rk (P_k) \|^2_F \label{eqn:app-7}\\
  & ={} \| P_k - P \|^2_F - 2 \eta \langle \nabla \rk (P_k), P_k - P \rangle +
  \eta^2 \| \nabla \rk (P_k) \|^2_F, \label{eqn:pf-lem4.1-1}
\end{align}
where \eqref{eqn:app-7} uses non-expansiveness of  projection. With convexity $\rk (P) - \rk (P_k) \geq \langle \nabla
\rk (P_k), P - P_k \rangle$,
\begin{align}
  \| P_{k + 1} - P \|^2_F \leq{} & \| P_k - P \|^2_F + 2 \eta (\rk (P) - \rk
  (P_k)) + \eta^2 \| \nabla \rk (P_k) \|^2_F \nonumber.
\end{align}
Re-arrangement yields $\rk (P_k) - \rk (P) \leq \tfrac{\eta}{2} \| \nabla \rk (P_k) \|^2 +
   \tfrac{1}{2 \eta} [\| P_k - P \|^2_F - \| P_{k + 1} - P \|^2_F]$.
Telescoping over $k$ and dropping the term $-\frac{1}{2 \eta} \| P_{K + 1} - P \|^2_F$ to obtain
\begin{equation} \label{eqn:pf-lem4.1-2}
\textstyle \sum_{k = 1}^K \rk (P_k) - \sum_{k = 1}^K \rk (P) \leq \tfrac{1}{2 \eta} \|
   P_1 - P \|^2_F + \tfrac{\eta}{2} \sum_{k = 1}^K \| \nabla \rk (P_k) \|^2_F.
\end{equation}

Using \Cref{lem:auxi-smooth}, we have $\textstyle \| \nabla \rk (P_k) \|_F^2 \leq 2\tau \rk (P_k)$. Plugging this bound into \eqref{eqn:pf-lem4.1-2} gives 
\begin{equation*}
\textstyle \sum_{k = 1}^K \rk (P_k) - \sum_{k = 1}^K \rk (P) \leq \tfrac{1}{2 \eta} \|
P_1 - P \|^2_F + \tau \eta \sum_{k = 1}^K \rk (P_k) .
\end{equation*}
Re-arrangement gives
\begin{equation} \label{eqn:pf-lem4.1-3}
\textstyle (1 - \tau \eta) \big[ \sum_{k = 1}^K \rk (P_k) - \sum_{k = 1}^K \rk (P) \big] \leq \tfrac{1}{2 \eta} \| P_1 - P \|^2_F +  \tau \eta \sum_{k = 1}^K \rk (P).
\end{equation}
For $\eta \leq \tfrac{1}{2 \tau}$, we may divide both sides of \eqref{eqn:pf-lem4.1-3} by $1 - \tau \eta$ and plug in the bound $\frac{1}{1 - \tau \eta} \leq 2$ to obtain
\begin{equation} \label{eqn:pf-lem4.1-4}
  \textstyle \sum_{k = 1}^K \rk (P_k) - \sum_{k = 1}^K \rk (P) 
\leq \tfrac{1}{\eta} \| P - P_1 \|^2_F + 2 \tau \eta \textstyle \sum_{k = 1}^K \rk (P),
\end{equation}
and this completes the proof.

\end{proof}

\subsection{Proof of Lemma \ref{lem:rx-surrogate-measure}}
Since $f(x) - f(x^\star) \geq 0$, applying \Cref{thm:online-to-conv-ratio} with $\varphi(x) = f(x) - f(x^\star)$ completes the proof.

\subsection{Proof of Proposition \ref{prop:rx-learnability}}

Denote $u_x (P) \assign f (x - P \nabla f(x))$. As a function of $P$, $u_x ( P ) = f (x - P \nabla f(x))$ is a composition of convex function $f (x)$ and affine transformation $x - P \nabla f(x)$. Hence $u_x$
is a convex function. The chain rule gives
\[ \nabla u_x ( P ) = \nabla f (x - P \nabla f(x)) = \nabla f (x - P \nabla f(x)) \nabla f(x)^{\top}.
\]
For any $P_1, P_2 \in \mathcal{P}$, we can successively deduce that
\begin{align}
\| \nabla u_x (P_1) - \nabla u_x (P_2) \|_F ={} & \| \nabla f (x -
P_1 \nabla f(x)) \nabla f(x)^{\top} - \nabla f (x - P_2 \nabla f(x)) \nabla f(x)^{\top} \|_F \nonumber\\
={} & \| [\nabla f (x - P_1 \nabla f(x)) - \nabla f (x - P_2 \nabla f(x))] \nabla f(x)^{\top} \|_F
\nonumber\\
\leq{} & \| \nabla f (x - P_1 \nabla f(x)) - \nabla f (x - P_2 \nabla f(x)) \| \cdummy \|  \nabla f(x) \| \label{eqn:app-1} \\
\leq{} & L \| (P_1 - P_2) \nabla f(x) \| \cdummy \|  \nabla f(x) \| \label{eqn:app-2} \\
\leq{} & L \|  \nabla f(x) \|^2 \| P_1 - P_2 \| \nonumber\\
\leq{} & L \|  \nabla f(x) \|^2 \| P_1 - P_2 \|_F, \nonumber
\end{align}

where \eqref{eqn:app-1} uses the submultiplicativity of Frobenius norm $\|AB\|_F \leq \|A\|_F \|B\|_F $; and \eqref{eqn:app-2} uses $L$-smoothness of $f(x)$. Hence $u_x$ is $L \|  \nabla f(x) \|^2$-smooth. Since the surrogate loss $r_x ( P) = \frac{u_x (P) - f(x^\star)}{f (x) - f(x^\star)}$ is a positive-scaled convex function $u_x$ with translation, and hence $r_x$ is also convex. 
Next, since $x \not \in \Xcal^\star$, the denominator of $r_x (P)$ must be positive, and hence $r_x (P) \geq 0$ for all $P$.
Lastly, since $r_x ( P) = \frac{u_x (P) - f(x^\star)}{f (x) - f(x^\star)}$ and $u_x$ is $L \| \nabla f (x) \|^2$-smooth, $r_x ( P)$ is also smooth with smoothness constant no greater than $2 L^2$:
\begin{equation*}
 \tfrac{L \| \nabla f (x) \|^2}{f
   (x) - f (x^{\star})} = \tfrac{L \| \nabla f (x) - \nabla f (x^{\star}) \|^2}{f
   (x) - f (x^{\star})} \leq 2 L^2,
\end{equation*}
where the last inequality invokes $L$-smoothness of $f(x)$. This completes the proof.

\subsection{Proof of Lemma \ref{lem:rx-learnability}}
For simplicity we denote $r_k(P) :=r_{x^k}(P)$.
By \Cref{prop:rx-learnability}, the surrogate losses $\{r_k\}$ are $2L^2$-smooth and non-negative. Then using \Cref{lem:smooth-online-learnability} with $\tau = 2L^2$, we get
\begin{equation} \label{eqn:pf-lem4.1-4}
  \textstyle \sum_{k = 1}^K \rk (P_k) - \sum_{k = 1}^K \rk (P) 
\leq \tfrac{1}{\eta} \| P - P_1 \|^2_F + 4 L^2 \eta \textstyle \sum_{k = 1}^K \rk (P),
\end{equation}
which proves \eqref{eqn:r-regret-bound}. Suppose further $\diam(\Pcal) \leq D$. Then 
\begin{align*}
r_x(P) &= \tfrac{f(x - P \nabla f(x)) - f(x^\star)}{f(x) - f(x^\star)} \\
&\leq \tfrac{f(x) - f(x^\star) - \langle \nabla f(x), (-P + \frac{L}{2} P^\top P) \nabla f(x) \rangle}{f(x) - f(x^\star)} \\
&\leq 1 + \| \tfrac{L}{2} P^\top P - P \| \tfrac{\| \nabla f(x) \|^2}{f(x) - f(x^\star)} \\
&\leq 1 + 2L (\tfrac{L}{2} D^2 + D) \\
&= (1 + L D)^2.
\end{align*}
Therefore, \eqref{eqn:pf-lem4.1-4} implies: for $\eta \leq \tfrac{1}{4L^2}$,
\begin{equation} \label{eqn:pf-lem4.1-5}
\textstyle \sum_{k = 1}^K \rk (P_k) - \sum_{k = 1}^K \rk (P) \leq \tfrac{1}{\eta} D^2 + 4 L^2 (1 + LD)^2 K \eta,
\end{equation}
in which the right-hand side is minimized (as a function of $\eta$) at $\eta = \tfrac{D}{2L (1 + LD) \sqrt{K}}$.
By taking the stepsize $\eta = \min\big\{ \tfrac{1}{4 L^2} , \tfrac{D}{2L (1 + LD) \sqrt{K}} \big\}$ and then minimizing over $P \in \mathcal{P}$, we conclude
\begin{equation*}
    \textstyle \sum_{k = 1}^K r_{x^k} (P_k) 
    \leq \displaystyle \min_{P \in \mathcal{P}} \textstyle \sum_{k = 1}^K r_{x^k} (P) + \max\big\{ 4 LD (1 + LD) \sqrt{K}, 8 L^2 D^2 \big\}.
\end{equation*}

and this completes the proof.

\subsection{Proof of Theorem \ref{thm:rx-trajectory-opt}}

By \Cref{lem:rx-learnability}, we have

\[\textstyle \frac{1}{K}\sum_{k = 1}^K r_{x^k} (P_k) \leq \frac{1}{K}\sum_{k = 1}^K r_{x^k} (P) +  \frac{\rho_K}{K}\]
for all $P \in \Pcal$. and plugging the relation into \Cref{lem:rx-surrogate-measure} completes the proof.

\subsection{Proof of Lemma \ref{lem:rx-hindsight}}
For any fixed $x \not \in \Xcal^\star$, the result $r_x (P_r^{\star}) \leq 1 - \frac{1}{\kappa^{\star}}$ is a direct consequence of \eqref{eqn:opt-contraction-factor}. 

\subsection{Proof of Corollary \ref{coro:rx-globalconv}}

Using \Cref{lem:rx-hindsight} and \Cref{thm:rx-trajectory-opt}, $\theta_K^\star \leq 1 - \frac{1}{\kappa^{\star}}$ and plugging the bound back into \Cref{thm:rx-trajectory-opt} completes the proof.

\subsection{Proof of Theorem \ref{thm:rx-globalconv}}
Combining \Cref{lem:rx-surrogate-measure} and \eqref{eqn:r-regret-bound} from \Cref{lem:rx-learnability} and using the relation $r_x(P^\star_r) \leq 1 - \tfrac{1}{\kappa^\star}$ from \Cref{lem:rx-hindsight}, we have, for $\eta \leq \tfrac{1}{4L^2}$, that
\begin{align}
f (x^{K + 1}) - f (x^{\star}) 
&\leq (f (x^1) - f (x^{\star}) ) \big(\textstyle \tfrac{1}{K} \sum_{k = 1}^K r_{x^k} (P^\star_r) + \tfrac{1}{\eta} \| P^\star_r - P_1 \|^2_F + 4 L^2 \eta \textstyle \sum_{k = 1}^K r_{x^k} (P^\star_r) \big)^K \nonumber \\
&\leq (f (x^1) - f (x^{\star}) ) \big(1 - \tfrac{1}{\kappa^\star} + \tfrac{1}{\eta} \| P^\star_r - P_1 \|^2_F + 4 L^2 K \eta \big)^K . \label{eqn:pf-thm4.2-1}
\end{align}
Take the stepsize $\eta = \min\big\{ \tfrac{1}{4 L^2} , \tfrac{\|P^{\star}_r - P_1\|_F}{2 L \sqrt{K}} \big\}$. The bound \eqref{eqn:pf-thm4.2-1} implies the desired result:
\begin{equation*}
    f (x^{K + 1}) - f (x^{\star}) \leq (f (x^1) - f (x^{\star})) \big( 1 - \tfrac{1}{\kappa^{\star}} + \max\big\{\tfrac{4 L \| P^{\star}_r - P_1 \|_F}{\sqrt{K}}, \tfrac{8L^2 \| P^{\star}_r - P_1 \|_F^2}{K} \big\} \big)^K.
\end{equation*}

\subsection{Proof of Proposition \ref{prop:rx-minimax}}
Given the optimization problem \eqref{eqn:minimax-sdp}, 
\begin{eqnarray*}
  \min_{P \in \mathcal{P}_+} \quad  \kappa \quad 
  \text{subject to} \quad \tfrac{1}{\kappa} I \preceq P^{1 / 2} A P^{1 / 2} \preceq  I,\label{eqn:quadratic-form}
\end{eqnarray*}
we can define $\tau = 1 / \kappa$ and reduce it to a standard semidefinite optimization problem (SDP)
\begin{equation}
  \max_{P \in \mathcal{P}_+} \quad \tau = \kappa^{-1} \quad
  \text{subject to} \quad A^{-1} \tau \preceq P \preceq A^{- 1}. \label{eqn:sdp}
\end{equation}
On the other hand, using $f(x) = \frac{1}{2}\langle x, A, x\rangle$, we can explicitly write 

\[r_x(P) = \frac{\frac{1}{2} \langle x, A (P A P - 2 P) A x \rangle}{\frac{1}{2}\langle x, A, x\rangle}\]
and $r_x$ is degree-zero homogeneous in $x$. Therefore, we can consider the following optimization problem
\begin{align*}
\min_{P \in \mathcal{P}_+} & \max_{\langle x, A x \rangle = 1}  ~~\langle x, A  (P A P - 2 P) A x \rangle, 
\end{align*}
which can be further re-written as
\begin{eqnarray*}
  \max_{P \in \mathcal{P}_+} \quad \lambda \quad
  \text{subject to} & 2 A^{1 / 2} P A^{1 / 2} - A^{1 / 2} P A P A^{1 / 2}
  \succeq \lambda I & 
\end{eqnarray*}
Next we do variable replacement by letting $M \assign A^{1 / 2} P A^{1 / 2}$ and $\Pcal'_+ = \{M = A^{1/2} P A^{1/2}: P \in \Pcal_+ \}$
and it suffices to show the equivalence between the following two problems.
\begin{equation} \tag{SDP} \label{eqn:app-sdp}
  \max_{M \in \mathcal{P}_+'} \quad  \tau  \quad
  \text{subject to} \qquad \tau I \preceq M \preceq I   \nonumber
\end{equation}
\begin{equation} \tag{Minimax} \label{eqn:app-minimax}
  \max_{M \in \mathcal{P}_+'} \quad \lambda \quad
  \text{subject to} \qquad  2 M - M^2 \succeq \lambda I \nonumber
\end{equation}
Given optimal solution $(M_1^{\star}, \tau^{\star})$ to \eqref{eqn:app-sdp}, we
have $\tau^{\star} I \preceq M^{\star}_1 \preceq I$ and let $M_1^{\star} = Q
\Lambda_1 Q^{\top}$. Plugging $M_1^{\star}$ into the constraint,
\[ 2 M_1^{\star} - (M_1^{\star})^2 = Q (2 \Lambda_1 - \Lambda_1^2) Q^{\top},
\]
which  corresponds to $\lambda = 2 \tau^{\star} - (\tau^{\star})^2$. On the
other hand, given optimal solution $(M_2^{\star}, \lambda^{\star})$ to \eqref{eqn:app-minimax}, $2 M_2^{\star} {- M_2^{\star}}^2 \succeq
\lambda^{\star} I$ and similarly we let $M_2^{\star} = Q \Lambda_2 Q^{\top}$. Then 
\[ Q (2 \Lambda_2 - \Lambda_2^2) Q^{\top} \succeq \lambda^{\star} I \]
and there exists some $q_j$ such that $2 \lambda_{2 j} - \lambda^2_{2 j} =
\lambda^{\star}$. It corresponds to $\tau = 2 \tau^2 - \tau =
\lambda^{\star}$. This establishes the equivalence between the two problems and completes the proof.

\subsection{Proof of Theorem \ref{thm:superlin-quadratic}}

Recall that by \eqref{eqn:pf-thm4.2-1} we have
\[ f (x^{K + 1}) - f (x^{\star}) \leq (f (x^1) - f (x^{\star})) (
   \tfrac{1}{K} \textstyle \sum_{k = 1}^K r_{x^k} (P^{\star}_r) + \tfrac{1}{\eta} \|
   P_r^{\star} - P_1 \|_F^2 + 4 L^2 \eta \textstyle \sum_{k = 1}^K r_{x^k} (P^{\star}_r)
   ) . \]
Using $r_{x^k} (P^{\star}_r) = 0, P^{\star}_r=A^{-1}$ and taking $\eta = 1 / (4 L^2)$, we get
\[ f (x^{K + 1}) - f (x^{\star}) \leq (f (x^1) - f (x^{\star})) (
   \tfrac{4 L^2 \| P_1 - A^{- 1} \|_F^2}{K} )^K . \]
This completes the proof.

\section{Function value ratio surrogate with optimal value lower bound} \label{sec:rxz-results}

This section analyzes the sub-optimal ratio surrogate loss $r^z_x (P)$ defined in \eqref{eqn:ratio-surrogate-loss-z}. Recall that
\begin{equation} \label{eqn:rx-subopt-loss}
r^z_x (P) \assign \tfrac{f (x - P \nabla f (x)) - z}{f (x) - z} = \tfrac{f (x^+) -
z}{f (x) - z},
\end{equation}
where $z < f(x^\star)$ is a lower bound for the optimal objective value. The challenging part of the analysis when $z < f(x^\star)$ is that the algorithm is only guaranteed to converge to some suboptimal solution whose suboptimality is determined by $f(x^\star) - z$, the accuracy of the lower bound. The analysis in this section is more involved than in \Cref{sec:ratio-surrogate}, and for clarity, we \textit{only present the global convergence result}.

\subsection{Surrogate loss}

\begin{lem}[Surrogate loss and measure] \label{lem:rzx-surrogate-measure}
For all $K \geq 1$, the online scaled gradient method satisfies
\begin{equation} \label{eqn:rzx-surrogate-measure}
f (x^{K + 1}) - z \leq ( f (x^1) - z) \big( \tfrac{1}{K} \textstyle \sum_{k = 1}^K \rxkz (P_k) \big)^K. 
\end{equation}
\end{lem}

\begin{prop}[Properties of $\rxkz$] \label{prop:rzx-learnability}
  Let $z < f(x^{\star})$ be a given lower bound.  Under \ref{A1} and \ref{A2}, for any fixed $x$, the surrogate loss $r^z_x (P)$ defined in \eqref{eqn:rx-subopt-loss} is convex, non-negative, and $2 L^2$-smooth as a function in $P$. In addition, the derivative of $r_x^z$ takes the form 
  \[\nabla r^z_x(P) = -\tfrac{\nabla f(x-P\nabla f(x)) \nabla f(x)^\top}{f(x) - z}.\]
\end{prop}

\subsection{Online learning algorithm}

\begin{lem}[Learnability] \label{lem:rzx-learnability}
Given \ref{A1}, \ref{A2}, and the ratio surrogate losses $\{\rxkz\}$,  online gradient descent
\begin{equation}\label{eqn:rxz-ogd}
	P_{k + 1} = \Pi_{\mathcal{P}} [P_k - \eta \nabla \rxkz (P_k)]
\end{equation}
with stepsize $\eta \leq 1/(4 L^2)$ generates a sequence of scaling matrices $\{P_k\}_{k \geq 2}$ such that 
\begin{equation} \label{eqn:r-regret-bound}
   \textstyle \sum_{k = 1}^K \rxkz (P_k) - \textstyle \sum_{k = 1}^K \rxkz (P) \leq \tfrac{1}{\eta} \| P - P_1 \|^2_F + 4 L^2 \eta \textstyle \sum_{k = 1}^K \rxkz (P) \quad \text{for any } P \in \mathcal{P}.
\end{equation}
\end{lem}

\subsection{Algorithm design and analysis}
In this section, we show how to obtain an $\Ocal(\kappa^\star \log^2 (1/\varepsilon))$ complexity through a double-loop algorithm. Since the double-loop algorithm deviates from our framework, only global convergence is established for brevity. We start by specifying the {\osgmrzx}, a subroutine that will be invoked in the inner loop.

We choose the optimality measure $\varphi$, the surrogate loss $\ell$, and the online learning algorithm $\mathcal{A}$ to be
\begin{equation*}
    \varphi(x) \assign f(x) - f(x^\star), \quad \ell_{x} (P) \assign \rxz (P), \quad \mathcal{A} \assign \text{online gradient descent in \eqref{eqn:rxz-ogd}},
\end{equation*}
and the monotone oracle $\Mcal$ is optional. \Cref{alg:osgm-rzx} presents \osgmrzx~without the monotone oracle.

\begin{algorithm}[h]
{\textbf{input $x^1, P_1 \in \Pcal, \eta >0, z < f(x^\star)$}}\\
\For{k =\rm{ 1, 2,...}}{
$x^{k + 1} = x^{k} - P_k \nabla f(x^k)$\\
$P_{k+1} = \Pi_{\Pcal} [P_k - \eta \nabla \rxkz (P_k) ]$

}
{\textbf{output} $x^\text{best}$ with minimum objective value}
\caption{Online scaled gradient method with lower bound ratio surrogate (\osgmrzx) \label{alg:osgm-rzx}}
\end{algorithm} 

\Cref{thm:unknown-fopt} characterizes the convergence behavior of {\osgmrzx}.
\begin{thm}[Global convergence with lower bound] \label{thm:unknown-fopt}
 Under \ref{A1} to \ref{A3}, \Cref{alg:osgm-rzx} (\osgmrzx) with $\eta
= \min\{ \tfrac{1}{4L^2}, \frac{\| P_r^{\star} - P_1 \|_F}{2 L \sqrt{K}} \}$ satisfies
  \begin{align}
    \min_{1 \leq k \leq K + 1} f (x^k) - f (x^{\star})
   & \leq{}  \tfrac{1}{2} (f (x^\star)
    - z) + (f (x^1) - f (x^{\star})) (1-\tfrac{1}{2\kappa^\star} + \tfrac{\rho_K}{K})^K,\nonumber
  \end{align}
where $\rho_K := \max\{4 L\sqrt{K} \| P^{\star}_r - P_1 \|_F, 8L \| P^{\star}_r - P_1 \|_F^2 \}$.
\end{thm}

\begin{lem}[Lower bound update] \label{lem:lbupdate}
  Under the same assumptions and parameter choice as \Cref{thm:unknown-fopt}
and denote \[z^+ = \frac{1}{2} \Big[\min_{1 \leq k \leq K + 1} f (x^k) + z \Big]. \] Then exactly
  one of the cases below happens:
  \begin{itemize}[leftmargin=15pt]
    \item $f (x^{K + 1}) - f (x^{\star})
    \leq  (f (x^1) - f (x^{\star})) (1-\tfrac{1}{2\kappa^\star} + \tfrac{\rho_K}{K})^K$, or 
    \item $f (x^{\star}) - z^+ \leq \tfrac{1}{2} (f (x^{\star}) - z)$ and $z^+
    \leq f (x^{\star})$.
  \end{itemize}
\end{lem}

\Cref{lem:lbupdate} suggests that the output of {\osgmrzx} either already satisfies the desirable convergence result, or the accuracy of the lower bound can be improved by a factor of 2. This motivates the idea of running {\osgmrzx} multiple times and outputting the best solution, as presented in \Cref{alg:double-loop}.
\begin{algorithm}[h]
{\textbf{input $x^1, P_1 \in \Pcal, \eta > 0, z^1 < f(x^\star)$}}\\
\For{t =\rm{ 1, 2,...}}{
$x^{t + 1} = \osgmrzx(x^t, P_1, \eta, z^t)$\\
$z^{t + 1} = \frac{1}{2}(f(x^{t + 1}) + z^t)$
}
{\textbf{output} $x^\text{best}$ with minimum objective value}
\caption{Online scaled gradient method with ratio surrogate and lower bound update\label{alg:osgd-lb}}
\label{alg:double-loop}
\end{algorithm} 

\Cref{thm:osgm-rxz} provides the final convergence result.
\begin{thm} \label{thm:osgm-rxz}
  Under the same assumptions and parameter choices as \tmtextbf{Theorem \ref{thm:unknown-fopt}},
  \tmtextbf{Algorithm \ref{alg:osgd-lb}} attains $f (x^{\tmop{best}}) - f (x^{\star}) \leq
  \varepsilon$ in at most $\mathcal{O} (\kappa^{\star} \log^2 (1 /
  \varepsilon))$ scaled gradient iterations.
\end{thm}

\subsection{Proof of Lemma \ref{lem:rzx-surrogate-measure}}

Since $f(x) - z > 0$, applying \Cref{thm:online-to-conv-ratio} with $\varphi(x) = f(x) - z$ completes the proof.

\subsection{Proof of Proposition \ref{prop:rzx-learnability}}
Since $z < f(x)$, both the numerator and the denominator of $\rxz$ are positive. Following the proof of \Cref{prop:rx-learnability}, we can show that $\rxz$ is convex. Since $u_x$ is $L\|\nabla f(x)\|^2$-smooth, $\rxz$ is $\tfrac{L \| \nabla f (x) \|^2}{f (x) - z}$-smooth, and
 
\begin{equation*}
\tfrac{L \| \nabla f (x) \|^2}{f (x) - z} < \tfrac{L \| \nabla f (x) \|^2}{f
   (x) - f (x^{\star})} = \tfrac{L \| \nabla f (x) - \nabla f (x^{\star}) \|^2}{f
   (x) - f (x^{\star})} \leq 2 L^2,
\end{equation*}
which completes the proof.

\subsection{Proof of Lemma \ref{lem:rzx-learnability}}
By \Cref{prop:rzx-learnability}, the surrogate losses $\{r_{x^k}^z\}$ are $2L^2$-smooth and non-negative. Applying \Cref{lem:smooth-online-learnability} with $\tau = 2L^2$ completes the proof.

\subsection{Proof of Theorem \ref{thm:unknown-fopt}}

Using \Cref{lem:rzx-surrogate-measure} and \Cref{lem:rzx-learnability} with $P = P_r^{\star}$,
\begin{align}
  f (x^{K + 1}) - z \leq{} & (f (x^1) - z) ( \tfrac{1}{K} \textstyle  \sum_{k = 1}^K
  r^z_{x^k} (P_k) )^K \nonumber\\
  \leq{} & (f (x^1) - z) ( \tfrac{1}{K} \textstyle \sum_{k = 1}^K r^z_{x^k}
  (P_r^{\star}) + \tfrac{1}{K} [ \tfrac{1}{\eta} \| P - P_r^{\star}
  \|_F^2 + 4 L^2 \eta \textstyle  \sum_{k = 1}^K r^z_{x^k} (P_r^{\star}) ] )^K
  \nonumber\\
  \leq{} & (f (x^1) - z) ( \tfrac{1}{K} \textstyle  \sum_{k = 1}^K r^z_{x^k}
  (P_r^{\star}) + \tfrac{1}{K} [ \tfrac{1}{\eta} \| P - P_r^{\star}
  \|_F^2 + 4 L^2 \eta ] )^K, \label{eqn:pf-c.1-1}
\end{align}

where \eqref{eqn:pf-c.1-1} uses $f (x - P_r^{\star} \nabla f (x)) \leq f (x)$ and that $r^z_x
(P_r^{\star}) = \tfrac{f (x - P_r^{\star} \nabla f (x)) - z}{f (x) - z} \leq
1$. 

Taking $\eta = \min \{ \tfrac{1}{4 L^2}, \tfrac{\| P_r^{\star} - P_1
\|_F}{2 L \sqrt{K}} \}$ gives
\begin{align*}
	f (x^{K + 1}) - z & \leq{} (f (x^1) - z) \big( \tfrac{1}{K} \textstyle \sum_{k = 1}^K
   r^z_{x^k} (P_r^{\star}) + \max \big\{ \tfrac{4 L \| P_r^{\star} - P_1
   \|_F}{\sqrt{K}}, \tfrac{8 L^2 \| P_r^{\star} - P_1 \|_F^2}{K}\big\}
   \big)^K \\ & ={} (f (x^1) - z) \big( \tfrac{1}{K} \textstyle \sum_{k = 1}^K
   r^z_{x^k} (P_r^{\star}) + \tfrac{\rho_K}{K} \big\}
   \big)^K.
\end{align*}
Next we analyze $\tfrac{1}{K} \sum_{k = 1}^K r^z_{x^k} (P_r^{\star})$, and
using
\[ f (x - P_r^{\star} \nabla f (x)) - f (x^{\star}) \leq ( 1 -
   \tfrac{1}{\kappa^{\star}} ) (f (x) - f (x^{\star})), \]
we deduce that
\begin{align}
  f (x - P_r^{\star} \nabla f (x)) - z ={} & f (x - P_r^{\star} \nabla f (x)) -
  f (x^{\star}) + f (x^{\star}) - z \nonumber\\
  \leq{} & ( 1 - \tfrac{1}{\kappa^{\star}} ) [f (x) - f (x^{\star})]
  + f (x^{\star}) - z \nonumber\\
  ={} & ( 1 - \tfrac{1}{\kappa^{\star}} ) [f (x) - z] - ( 1 -
  \tfrac{1}{\kappa^{\star}} ) [f (x^{\star}) - z] + f (x^{\star}) - z
  \nonumber\\
  ={} & ( 1 - \tfrac{1}{\kappa^{\star}} ) [f (x) - z] +
  \tfrac{1}{\kappa^{\star}} [f (x^{\star}) - z] \nonumber
\end{align}

Dividing both sides by $f (x) - z$ gives
\[ \tfrac{f (x - P_r^{\star} \nabla f (x)) - z}{f (x) - z} = ( 1 -
   \tfrac{1}{\kappa^{\star}} ) + \tfrac{1}{\kappa^{\star}} \tfrac{f
   (x^{\star}) - z}{f (x) - z} = 1 - \tfrac{1}{\kappa^{\star}} \tfrac{f (x) -
   f (x^{\star})}{f (x) - z} . \]
Hence $\tfrac{1}{K} \sum_{k = 1}^K r^z_{x^k} (P_r^{\star}) \leq 1 -
\frac{1}{\kappa^{\star}} ( \tfrac{1}{K} \sum_{k = 1}^K \frac{f (x^k) - f
(x^{\star})}{f (x^k) - z} )$. Now, we do case analysis

\paragraph{Case 1.} Suppose $\frac{f (x^k) - f (x^{\star})}{f (x^k) - z} \geq
\frac{1}{2}$ for all $1 \leq k \leq K$, then $\tfrac{1}{K} \sum_{k = 1}^K
r^z_{x^k} (P_r^{\star}) \leq 1 - \frac{1}{2 \kappa^{\star}}$ and
\begin{align}
  \min_{1 \leq k \leq K + 1} f (x^k) - f (x^{\star}) \leq{} & f (x^{K + 1}) - f
  (x^{\star}) \nonumber\\
  \leq{} & (f (x^1) - f (x^{\star})) \big( 1 - \tfrac{1}{2 \kappa^{\star}} +
  \tfrac{\rho_K}{K} \big)^K . \nonumber
\end{align}
\paragraph{Case 2.} Otherwise, there is some $1 \leq j \leq K$ such
that $\frac{f (x^j) - f (x^{\star})}{f (x^j) - z} \leq \frac{1}{2}$, a
re-arrangement gives $2 f (x^j) - 2 f (x^{\star}) \leq f (x^j) - z$ and
\[ \min_{1 \leq k \leq K + 1} f (x^k) - f (x^{\star}) \leq f (x^j) - f
   (x^{\star}) \leq \tfrac{1}{2} (f (x^j) - z) . \]
Putting the two cases together, we complete the proof.

\subsection{Proof of Lemma \ref{lem:lbupdate}}

The argument is the same as in \tmtextbf{Theorem \ref{thm:unknown-fopt}}. In \tmtextbf{Case 1}, we
get the first convergence result. Otherwise, we know that there exists some $1
\leq j \leq K$ such that $\frac{f (x^j) - f (x^{\star})}{f (x^j) - z} \leq
\frac{1}{2}$ and since $\min_{1 \leq k \leq K} f (x^k) \leq f (x^j)$, we have
\[ \tfrac{\min_{1 \leq k \leq K} f (x^k) - f (x^{\star})}{\min_{1 \leq k \leq
   K} f (x^k) - z} \leq \tfrac{f (x^j) - f (x^{\star})}{f (x^j) - z} \leq
   \tfrac{1}{2} . \]
Rearranging the relation, we have $z^+ = \frac{1}{2} [\min_{1 \leq k \leq K} f
(x^k) + z] \leq f (x^{\star})$ and
\begin{align}
  f (x^{\star}) - \tfrac{1}{2} \big[\min_{1 \leq k \leq K} f (x^k) + z \big] ={} &
  \tfrac{1}{2} [f (x^{\star}) - \min_{1 \leq k \leq K} f (x^k)] + \tfrac{1}{2}
  [f (x^{\star}) - z] \nonumber\\
  \leq{} & \tfrac{1}{2} [f (x^{\star}) - z] . \nonumber
\end{align}
This completes the proof.

\subsection{Proof of Theorem \ref{thm:osgm-rxz}}

Denote $x^{t+1}$ as the output of {\osgmrzx} in iteration $t$ of \tmtextbf{Algorithm \ref{alg:osgd-lb}}.
If we fall into \tmtextbf{Case 1} in \tmtextbf{Lemma \ref{lem:lbupdate}} after running {\osgmrzx} for $K$ iterations, then 
\begin{equation} \label{eqn:pf-thm-osgm-1}
   f (x^{t+1}) - f (x^{\star}) \leq
   (f (x^t) - f (x^{\star})) \big( 1 - \tfrac{1}{2 \kappa^{\star}} +
  \tfrac{\rho_K}{K} \big)^K
\end{equation}
since $x^t$ is the initial point of {\osgmrzx} in iteration $t$.
Using the fact that $z^{1}$ is a lower bound for $f(x^{\star})$, algebraic manipulation shows that the right-hand side of \eqref{eqn:pf-thm-osgm-1} is less than $\varepsilon$ whenever
\begin{equation*}
   K \geq \tfrac{128 \|P^{\star}_r - P_1\|_F^2 L^2}{\log^2(1-\frac{1}{2\kappa^{\star}})} + 2 \kappa^{\star} \log\big( \tfrac{f (x^1) - z^{1}}{\varepsilon} \big) =: K_0.
\end{equation*}
We claim that if we run {\osgmrzx} for $K_0$ iterations at each iteration $t$ in \tmtextbf{Algorithm \ref{alg:osgd-lb}} and run \tmtextbf{Algorithm \ref{alg:osgd-lb}} for $T := \tfrac{1}{\log 2} \log\big( \tfrac{4(f(x^1) - z^1)}{\varepsilon} \big)$ iterations, then we have $f(x^{\text{best}}) - f(x^{\star}) \leq \varepsilon$ where $x^{\text{best}}$ is the point in $\{x^t: t = 1, \ldots, T+1\}$ with the smallest function value.
Hence, \tmtextbf{Algorithm \ref{alg:osgd-lb}} takes at most $K_0 T = \mathcal{O}(\kappa^{\star} \log^2(1/\varepsilon))$ scaled gradient iterations.\\

We will show that at least one of the iterates in $\{x^t: t = 1, \ldots, T+1\}$ from our algorithm satisfies $f(x^t) - f(x^{\star}) < \varepsilon$.
If we fall into \tmtextbf{Case 1} in \tmtextbf{Lemma \ref{lem:lbupdate}} for some iteration $t$, then we have $f(x^{t+1}) - f(x^{\star}) \leq \varepsilon$ by \eqref{eqn:pf-thm-osgm-1}.
Otherwise, we fall into \tmtextbf{Case 2} in \tmtextbf{Lemma \ref{lem:lbupdate}} for all $t \leq T$. 
In this case, we halve the distance between $z^t$ and $f(x^{\star})$ after every outer iteration, so that after $T := \big\lceil \tfrac{1}{\log 2} \log\big( \tfrac{4(f(x^1) - z^1)}{\varepsilon} \big) \big\rceil$ iterations, we have
\begin{equation*}
   |z^T - f(x^{\star})| \leq \left(\tfrac{1}{2} \right)^{T-1} (f(x^{\star})-z^1) \leq \left(\tfrac{1}{2} \right)^{T-1} (f(x^{1})-z^1) \leq \tfrac{\varepsilon}{2}.
\end{equation*}
Since $z^{T+1} = \frac{1}{2}(f(x^{t+1}) + z^T)$ and we fall into \tmtextbf{Case 2} at iteration $T$, we have
\begin{equation*}
|z^{T+1} - f(x^{\star})| = 
\big| \tfrac{f(x^{T+1}) + z^T}{2} - f(x^{\star}) \big| \leq \tfrac{1}{2} | z^{T} - f(x^{\star})| \leq \tfrac{\varepsilon}{4}.
\end{equation*}
Rearranging the relation, we have $f(x^{T+1}) \leq f(x^{\star}) + \frac{1}{2} \varepsilon + (f(x^{\star}) - z^T) \leq f(x^{\star}) + \varepsilon$. This completes the proof.

\section{Proof of results in Section \ref{sec:gnorm-surrogate}}

\subsection{Proof Lemma \ref{lem:gx-surrogate-measure}}

Given monotone oracle $\mathcal{M}$ with respect to gradient norm and by
definition of $g_x$,
\[ \| \nabla f (x^{k + 1}) \| = \| \nabla f (\mathcal{M} (x^{k})) \|
   \leq \| \nabla f (x^{k} - P_k \nabla f(x^k)) \| = g_{x^k} (P_k) \| \nabla f (x^k) \| .
\]
Hence, through the same argument as in \Cref{thm:online-to-conv-ratio}, we deduce that
\[ \tfrac{\| \nabla f (x^{K + 1}) \|}{\| \nabla f (x^1) \|} = \textstyle \prod_{k = 1}^K
   \tfrac{\| \nabla f (x^{k + 1}) \|}{\| \nabla f (x^k) \|} \leq (
   \tfrac{1}{K} \sum_{k = 1}^K \tfrac{\| \nabla f (x^{k + 1}) \|}{\| \nabla f
   (x^k) \|} )^K \leq ( \tfrac{1}{K} \sum_{k = 1}^K g_{x^k} (P_k)
   )^K \]
and this completes the proof.

\subsection{Proof of Proposition \ref{prop:gx-learnability}}

Lipschitzness of $g_x$ is straight-forward:
\begin{align}
  | g_x (P_1) - g_x (P_2) | ={} & \Big| \tfrac{\| \nabla f (x - P_1 \nabla f
  (x)) \|}{\| \nabla f (x) \|} - \tfrac{\| \nabla f (x - P_2 \nabla f (x))
  \|}{\| \nabla f (x) \|} \Big| \nonumber\\
  \leq{} & \tfrac{\| \nabla f (x - P_1 \nabla f (x)) - \nabla f (x - P_2 \nabla
  f (x)) \|}{\| \nabla f (x) \|} \label{eqn:app-12} \\
  \leq{} & \tfrac{L \| P_1 - P_2 \| \cdummy \| \nabla f (x) \|}{\| \nabla f (x)
  \|} \label{eqn:app-13} \\
  ={} & L \| P_1 - P_2 \| \leq L \| P_1 - P_2\|_F \, \nonumber
\end{align}
where \eqref{eqn:app-12} uses the triangle inequality $|\|a\| - \|b\|| \leq \|a - b\|$ and \eqref{eqn:app-13} uses $L$-smoothness of $f$. Next consider $| g_x (P) - \hat{g}_x (P) |$ and we deduce that
\begin{align}
  | g_x (P) - \hat{g}_x (P) | ={} & \Big| \Big\| \tfrac{\nabla f (x)}{\|
  \nabla f (x) \|} - \textstyle \int_0^1 \nabla^2 f (x - t P \nabla f (x)) P
  \tfrac{\nabla f (x)}{\| \nabla f (x) \|} ~\mathd t \Big\| - \Big\| \tfrac{\nabla f
  (x)}{\| \nabla f (x) \|} - \nabla^2 f (x) P \tfrac{\nabla f (x)}{\| \nabla f
  (x) \|} \Big\| \Big| \nonumber\\
  \leq{} & \| \textstyle \int_0^1 \nabla^2 f (x - t P \nabla f (x)) P \tfrac{\nabla f
  (x)}{\| \nabla f (x) \|} - \nabla^2 f (x) P \tfrac{\nabla f (x)}{\| \nabla f
  (x) \|} ~\mathd t \| \label{eqn:app-14} \\
  \leq{} & \textstyle \int_0^1 \| \nabla^2 f (x - t P \nabla f (x)) - \nabla^2 f (x) \|
  ~\mathd t \cdot \big (\tfrac{\| P \| \| \nabla f (x) \|}{\| \nabla f (x) \|} \big) \label{eqn:app-15} \\
  \leq{} & H \textstyle \int_0^1  \| P \nabla f (x) \| t ~\mathd t \cdot \| P \| \label{eqn:app-16} \\
  \leq {} & \tfrac{1}{2}H \| P \|^2 \| \nabla f (x) \|  \label{eqn:pf-prop-4.1-0}
\end{align}

where \eqref{eqn:app-14} again uses  $|\|a\| - \|b\|| \leq \|a - b\|$ and  \eqref{eqn:app-15} uses the Lispchitzness of the Hessian; \eqref{eqn:app-16} uses $\|A B \| \leq \|A \|\|B\|$. Convexity of $\hat{g}_x$ is straight-forward since $\hat{g}_x$ is a composition of linear function (in $P$) with norm $\| \cdummy \|$.
To show $L$-Lipschitzness of $\hat{g}_x$, we have
\begin{align}
  | \hat{g}_x (P_1) - \hat{g}_x (P_2) | ={} & \Big| \Big\| \tfrac{\nabla f
  (x)}{\| \nabla f (x) \|} - \nabla^2 f (x) P_1 \tfrac{\nabla f (x)}{\| \nabla
  f (x) \|} \Big\| - \Big\| \tfrac{\nabla f (x)}{\| \nabla f (x) \|} -
  \nabla^2 f (x) P_2 \tfrac{\nabla f (x)}{\| \nabla f (x) \|} \Big\| \Big|
  \nonumber\\
  \leq{} & \tfrac{1}{\| \nabla f (x) \|} \| \nabla^2 f (x) (P_1 - P_2) \nabla f
  (x) \| \label{eqn:pf-prop-4.1-1} \\
  \leq{} & L \| P_1 - P_2 \| \label{eqn:pf-prop-4.1-2},
\end{align}
where \eqref{eqn:pf-prop-4.1-1} again uses $|\|a\| - \|b\|| \leq \|a - b\|$ and \eqref{eqn:pf-prop-4.1-2} uses $\|\nabla^2 f(x)\| \leq L$. Last, we combine the convex subgradient lower bound of $\hat{g}_x$ with the approximation
\begin{align}
  g_x (P_1) \geq{} & \hat{g}_x (P_1) - \tfrac{1}{2} H \| P_1 \|^2 \| \nabla f (x)
  \| \label{eqn:pf-prop-4.1-3} \\
  \geq{} & \hat{g}_x (P_2) + \langle \hat{g}_x' (P_2), P_1 - P_2 \rangle -
  \tfrac{1}{2} H \| P_1 \|^2 \| \nabla f (x) \| \label{eqn:pf-prop-4.1-4} \\
  \geq{} & g_x (P_2) + \langle \hat{g}_x' (P_2), P_1 - P_2 \rangle - \tfrac{1}{2} H [\| P_1 \|^2 + \| P_2 \|^2]
  \| \nabla f (x) \| \label{eqn:pf-prop-4.1-5}\\
  \geq{} & g_x (P_2) + \langle \hat{g}_x' (P_2), P_1 - P_2 \rangle - H D^2
  \| \nabla f (x) \| \nonumber,
\end{align}
where \eqref{eqn:pf-prop-4.1-3} uses \eqref{eqn:pf-prop-4.1-0}, \eqref{eqn:pf-prop-4.1-4} uses convexity of $\hat{g}_x$ and \eqref{eqn:pf-prop-4.1-5} again applies \eqref{eqn:pf-prop-4.1-0}. This completes the proof.

\subsection{Proof of Lemma \ref{lem:gx-learnability}}

Denote $g_k (P_k) \assign g_{x^k} (P_k)$. For any $P \in
\mathcal{P}$, we have
\begin{align}
  \| P_{k + 1} - P \|^2_F & = \| \Pi_{\mathcal{P}} [P_k - \eta {g}_k'
  (P_k)] \|^2_F \nonumber\\
  & \leq \| P_k - P - \eta {g}_k' (P_k) \|^2_F \nonumber\\
  & = \| P_k - P \|^2_F - 2 \eta \langle {g}_k' (P_k), P_k - P \rangle +
  \eta^2 \| {g}_k' (P_k) \|^2_F \nonumber\\
  & \leq \| P_k - P \|^2_F - 2 \eta [g_k (P_k) - g_k (P)] + \eta^2 L^2 + \eta
  H D^2 \| \nabla f (x^k) \|, \label{eqn:pf-lem-5.2-1}
\end{align}
where \eqref{eqn:pf-lem-5.2-1} invokes \Cref{prop:hx-learnability}. 
Dividing both sides by $\eta$ and re-arranging the terms,
\begin{align}
  g_k (P_k) - g_k (P) \leq{} & \tfrac{\| P_k - P \|^2_F}{2 \eta} - \tfrac{\|
  P_{k + 1} - P \|^2_F}{2 \eta} + \tfrac{\eta}{2} L^2 + \tfrac{H D^2}{2} \|
  \nabla f (x^k) \| \nonumber\\
  \leq{} & \tfrac{\| P_k - P \|^2_F}{2 \eta} - \tfrac{\| P_{k + 1} - P \|^2_F}{2
  \eta} + \tfrac{\eta}{2} L^2 + \tfrac{H D^2}{2} \| \nabla f (x^1) \|, \label{eqn:pf-lem-5.2-2}
\end{align}

where \eqref{eqn:pf-lem-5.2-2} uses we use the fact that $\| \nabla f (x^k) \| \leq \| \nabla f (x^1) \|$. Summing both sides from $k = 1, \ldots, K$, we get the desired result:
\begin{align}
  \textstyle \sum_{k = 1}^K g_k (P_k) - \sum_{k = 1}^K g_k (P) \leq{} & \tfrac{\| P_1 - P
  \|^2_F}{2 \eta} + \tfrac{\eta}{2} L^2 K + \tfrac{H D^2}{2} \| \nabla f
  (x^1) \| K.
\end{align}
Finally, using the bound $\| P_1 - P
  \|^2_F \leq 4 D^2$ and plugging in the stepsize $\eta = \tfrac{2D}{L \sqrt{K}}$ yield \eqref{eqn:gx-regret}.

\subsection{Proof of Theorem \ref{thm:gx-trajectory-opt}}

By \Cref{lem:gx-learnability}, we have
\[\textstyle \frac{1}{K}\sum_{k = 1}^K g_{x^k} (P_k) \leq \displaystyle \min_{P \in \mathcal{P}} \tfrac{1}{K} \textstyle \sum_{k = 1}^K g_{x^k} (P) +  \frac{\rho_K}{K}\]
Plugging the relation into \Cref{lem:gx-surrogate-measure} completes the proof.

\subsection{Proof of Lemma \ref{lem:gx-hindsight}}

The first relation follows from
\begin{equation} \label{eqn:pf-lem5.3-1}
\| \nabla f (x - P_g^\star \nabla f (x)) \| \leq [ \textstyle\int_0^1 ~\| I - \nabla^2 f (x - t P_g^\star \nabla f (x)) P_g^\star \| \mathd
  t ] \cdot \| \nabla f (x) \| \leq \omega^{\star} \| \nabla f (x) \| = (1 - \tfrac{1}{\lambda^{\star}}) \| \nabla f (x) \| .
\end{equation}
The fact that $\frac{\mu}{L} I \preceq L^{- 1} \nabla^2 f (x) \preceq I$ for all $x$ implies $\| I - L^{- 1} \nabla^2 f (x) \| \leq 1 - \tfrac{\mu}{L}$ for all $x$. Hence, by taking $P = \frac{1}{L} I \in \mathcal{P}$, we conclude
\begin{equation*}
    \omega^{\star} \assign \min_{P \in \mathcal{P}} \max_x  \| I - \nabla^2 f(x) P \| \leq \max_x  \| I - L^{-1} \nabla^2 f
   (x) \| \leq 1 - \tfrac{\mu}{L}.
\end{equation*}
The desired inequality $\lambda^\star \leq \tfrac{L}{\mu}$ immediately follows from the definition $\lambda^\star = \frac{1}{1- \omega^\star}$.
Finally, rearranging \eqref{eqn:pf-lem5.3-1} gives the desired bound on gradient norm surrogate loss:
\begin{equation*}
    g_x(P_g^\star) = \tfrac{\| \nabla f (x - P_g^\star \nabla f (x)) \|}{\| \nabla f (x) \|} \leq 1 - \tfrac{1}{\lambda^{\star}}.
\end{equation*}

\begin{rem}
We can link $\lambda^\star$ and $\kappa^\star$ through two relations below:
\begin{align}
\| \nabla f (x - P^{\star}_g \nabla f (x)) \| \leq{} & (1 - \tfrac{1}{\lambda^{\star}}) \| \nabla f (x) \|, \nonumber \\
\| \nabla f (x - P_r^{\star} \nabla f (x)) \|_{P_r^{\star}} \leq{} & (
    1 - \tfrac{1}{\kappa^{\star}} ) \| \nabla f (x) \|_{P_r^{\star}}. \label{eqn:Pr-contraction}
\end{align}
The second relation \eqref{eqn:Pr-contraction} holds by simple algebraic derivation: since $\tfrac{1}{\kappa^{\star}} I \preceq
(P_r^{\star})^{1 / 2} \nabla^2 f (x) (P_r^{\star})^{1 / 2} \preceq I$ for all $x$, we  deduce that
\begin{align}
  & \| \nabla f (x - P_r^{\star} \nabla f (x)) \|_{P_r^{\star}}^2 \nonumber\\
  ={} & \| \nabla f (x) - \textstyle \int_0^1 \nabla^2 f (x - t P_r^{\star} \nabla f
  (x)) P_r^{\star} \nabla f (x) \mathd t \|_{P_r^{\star}}^2 \nonumber\\
  ={} & \| \textstyle \int_0^1 (I - \nabla^2 f (x - t P_r^{\star} \nabla f (x))
  P_r^{\star}) \nabla f (x) ~\mathd t \|_{P_r^{\star}}^2 \nonumber\\
  ={} & \langle \textstyle \int_0^1 (I - \nabla^2 f (x - t P_r^{\star} \nabla f (x))
  P_r^{\star}) \nabla f (x) ~\mathd t, P_r^{\star}  \textstyle \int_0^1 (I - \nabla^2 f (x - t
  P_r^{\star} \nabla f (x)) P_r^{\star}) \nabla f (x) ~\mathd t \rangle
  \nonumber\\
  ={} & \langle (P_r^{\star})^{1 / 2} \nabla f (x), ( \textstyle \int_0^1 (I - M_t)
  ~\mathd t )^2 (P_r^{\star})^{1 / 2} \nabla f (x) \rangle, \nonumber
\end{align}
where $M_t \assign (P_r^{\star})^{1 / 2} \nabla^2 f (x - t P^\star \nabla f(x) ) (P_r^{\star})^{1 / 2}$. Using the fact that $\tfrac{1}{\kappa^{\star}} \preceq M_t \preceq I$, we have $\int_0^1 (I - M)^2 \mathd t \preceq
( 1 - \tfrac{1}{\kappa^{\star}} )^2 I$ and hence
\[ \| \nabla f (x - P_r^{\star} \nabla f (x)) \|_{P_r^{\star}}^2 \leq ( 1 -
   \tfrac{1}{\kappa^{\star}} )^2 \| \nabla f (x) \|_{P_r^{\star}}^2 . \]
Taking square root on both sides gives the desired relation. However, since evaluating $\|\cdot \|_{P_r^{\star}}$ requires knowledge of $P_r^{\star}$, we have to define auxiliary quantity $\lambda^\star$ and $P_g^{\star}$.
\end{rem}

\subsection{Proof of Corollary \ref{coro:gx-globalconv}}
With convergence results from standard gradient descent, it takes $\Ocal(\kappa \log (H D^2\lambda^\star))$ iterations to output $\hat{x}$ such that $\|\nabla f(\hat{x})\| \leq \frac{1}{H D^2\lambda^\star}$. Next let $x^1 = \hat{x}$. Using \Cref{lem:gx-hindsight} and \Cref{thm:gx-trajectory-opt}, $\theta_K^{\star} \leq
1 - \tfrac{1}{\lambda^{\star}}$ and
\begin{align}
  \| \nabla f (x^{K + 1}) \| \leq{} & \| \nabla f (x^1) \| ( 1 -
  \tfrac{1}{\lambda^{\star}} + \tfrac{2 D L}{\sqrt{K}} + \tfrac{H D^2}{2} \|
  \nabla f (x^1) \| )^K \nonumber\\
  \leq{} & \| \nabla f (x^1) \| ( 1 - \tfrac{1}{\lambda^{\star}} + \tfrac{2
  D L}{\sqrt{K}} + \tfrac{1}{2 \lambda^{\star}} )^K \label{eqn:pf-coro-5.1-1}\\
  \leq{} & \| \nabla f (x^1) \| ( 1 - \tfrac{1}{2 \lambda^{\star}} +
  \tfrac{2 D L}{\sqrt{K}} )^K, \nonumber
\end{align}

where \eqref{eqn:pf-coro-5.1-1} uses the assumption that $\| \nabla f (x^1) \| \leq \tfrac{1}{H D^2 \lambda^{\star}}$. This completes the proof.

% \subsection{Improved analysis using adaptive subgradient methods}

\section{Proof of results in Section \ref{sec:hypergrad-surrogate}}

\subsection{Proof of Lemma \ref{lem:hx-surrogate-measure}}

For convenience we denote $x^{k+1/2} := x^k - P_k \nabla f(x^k)$. By definition of the monotone oracle, we always have 
\[f (x^{k + 1}) = f (\mathcal{M} (x^{k})) \leq \min \{ f
(x^k), f (x^{k + 1 / 2}) \}.\]
\paragraph{Proof of relation \eqref{eqn:hx-strong-cvx}.}
Suppose $\mu \neq 0$. By definition of $h_x (P)$, we can write
\begin{align}
  f (x^{k + 1 / 2}) - f (x^{\star}) & =  f (x^k) - f (x^{\star}) + h_{x^k}
  (P_k) \| \nabla f (x^k) \|^2 \nonumber \\
  & =  (f (x^k) - f (x^{\star})) \Big[ 1 + \tfrac{h_{x^k} (P_k) \| \nabla f
  (x^k) \|^2}{f (x^k) - f (x^{\star})} \Big] . \label{eqn:pf-lem6.1-1}
\end{align}
Since $f (x^{k + 1}) = f (\mathcal{M} (x^{k})) \leq \min \{ f
(x^k), f (x^{k + 1 / 2}) \}$, we have
\begin{align}
  f (x^{k + 1}) - f (x^{\star})
  \leq{} & \min\{f(x^k) - f(x^{\star}), f(x^{k+1/2}) - f(x^{\star}) \} \nonumber\\
  ={} & (f (x^k) - f (x^{\star})) \Big[ 1 + \min \Big\{ \tfrac{h_{x^k}
  (P_k) \| \nabla f (x^k) \|^2}{f (x^k) - f (x^{\star})}, 0 \Big\} \Big] \label{eqn:pf-lem-6.1-1},
\end{align}
where \eqref{eqn:pf-lem-6.1-1} uses \eqref{eqn:pf-lem6.1-1}. We successively deduce that
\begin{align}
  \tfrac{f (x^{K + 1}) - f (x^{\star})}{f (x^1) - f (x^{\star})} ={} & \textstyle \prod_{k
  = 1}^K \tfrac{f (x^{k + 1}) - f (x^{\star})}{f (x^k) - f (x^{\star})}
  \nonumber\\
  \leq{} & \big( \tfrac{1}{K} \textstyle \sum_{k=1}^K \tfrac{f(x^{k+1}) - f(x^{\star})}{f(x^{k}) - f(x^{\star})} \big)^K \nonumber \\
  \leq{} & \Big( 1 + \tfrac{1}{K} \textstyle\sum_{k = 1}^K \min \Big\{ \tfrac{h_{x^k}
  (P_k) \| \nabla f (x^k) \|^2}{f (x^k) - f (x^{\star})}, 0 \Big\} \Big)^K
  \nonumber\\
  \leq{} & \big( 1 + \tfrac{2 \mu}{K} \textstyle\sum_{k = 1}^K \min \{ h_{x^k} (P_k), 0
  \} \big)^K, \label{eqn:pf-lem-6.1-2}
\end{align}

where \eqref{eqn:pf-lem-6.1-2} is by $\tfrac{1}{2 \mu} \| \nabla f (x^k) \|^2 \geq
f (x^k) - f (x^{\star})$ and $\min \{ h_{x^k} (P_k), 0
  \} \leq 0$: 
\[ \min \Big\{ \tfrac{h_{x^k} (P_k) \| \nabla f (x^k) \|^2}{f (x^k) - f
   (x^{\star})}, 0 \Big\} = \tfrac{\| \nabla f (x^k) \|^2}{f (x^k) - f
   (x^{\star})} \cdot \min \{ h_{x^k} (P_k), 0 \} \leq 2 \mu \min \{ h_{x^k}
   (P_k), 0 \} . \]
By concavity of $\min \{ \cdot, 0 \}$, we have
\[ 1 + \tfrac{2 \mu}{K} \textstyle \sum_{k = 1}^K \min \{ h_{x^k} (P_k), 0 \} \leq 1 + 2
   \mu \min \{ \tfrac{1}{K} \sum_{k = 1}^K h_{x^k} (P_k), 0 \} \]
and using the identity $\max \{\cdot , 0\} = -\min \{-(\cdot) , 0\}$ completes the proof.

\paragraph{Proof of relation \eqref{eqn:hx-cvx-gnorm}.}
Again by definition of $h_x(P)$, $f (x^{k + 1 / 2}) - f (x^k) =
h_{x^k} (P_k) \| \nabla f (x^k) \|^2$ and
\begin{align}
  f (x^{k + 1}) - f (x^k) \leq{} & \min \big\{ f (x^{k + 1 / 2}) - f (x^k), f (x^{k}) - f (x^k) \big\}
  ={} \min \{ h_{x^k} (P_k), 0 \} \| \nabla f (x^k) \|^2 . \nonumber
\end{align}
Summing the inequality from $k=1$ to $K$, we have
\[ f (x^{K + 1}) - f (x^1) \leq
  \textstyle \sum_{k = 1}^K \min \{ h_{x^k} (P_k), 0 \} \| \nabla f (x^k) \|^2. \]
Re-arrangement gives
\begin{align}
  & \textstyle ( \sum_{k = 1}^K \max \{ - h_{x^k} (P_k), 0 \} ) \cdot \displaystyle \min_{1 \leq
  k \leq{} K}  \| \nabla f (x^k) \|^2 \nonumber\\
  \leq{} & \textstyle \sum_{k = 1}^K \max \{ - h_{x^k} (P_k), 0 \} \| \nabla f (x^k) \|^2
  \nonumber\\
  \leq{} & f (x^1) - f (x^{K + 1}) \nonumber\\
  \leq{} & f (x^1) - f (x^{\star}) \nonumber
\end{align}
Last, using convexity of $\max \{ \cdot, 0 \}$, 
\begin{align}
  \min_{1 \leq k \leq K}  \| \nabla f (x^k) \|^2 \leq{} & \tfrac{f (x^1) - f
  (x^{\star})}{K} \tfrac{1}{\tfrac{1}{K} \sum_{k = 1}^K \max \{ - h_{x^k}
  (P_k), 0 \}} \nonumber\\
  \leq{} & \tfrac{f (x^1) - f (x^{\star})}{K} \tfrac{1}{\max \{
  \frac{1}{K} \sum_{k = 1}^K - h_{x^k} (P_k), 0 \}} \nonumber
\end{align}
and this completes the proof.

\paragraph{Proof of relation \eqref{eqn:hx-cvx-fval}.}
Take $x^\star$ to be the equilibrium of the inner problem $\max_{x \in \Lcal_{f(x^1)}} \min_{x^{\star} \in
  \mathcal{X}^{\star}} \| x - x^{\star} \|$, we deduce that
\begin{align}
  f (x^{k + 1}) - f (x^{\star}) \leq{} & f (x^k) - f (x^{\star}) + \min \{
  h_{x^k} (P_k), 0 \} \| \nabla f (x^k) \|^2 \nonumber\\
  ={} & f (x^k) - f (x^{\star}) + \min \{ h_{x^k} (P_k), 0 \}  \tfrac{\| \nabla
  f (x^k) \|^2 \| x^k - x^{\star} \|^2}{(f (x^k) - f (x^{\star}))^2} \tfrac{[f
  (x^k) - f (x^{\star})]^2}{\| x^k - x^{\star} \|^2} \nonumber\\
  \leq{} & f (x^k) - f (x^{\star}) + \min \{ h_{x^k} (P_k), 0 \}  \centering\tfrac{[f
  (x^k) - f (x^{\star})]^2}{\| x^k - x^{\star} \|^2}, \label{eqn:pf-lem-6.1-3}
\end{align}

where the last inequality uses $f (x^k) - f (x^{\star}) \leq \| \nabla f (x^k)
\| \| x^k - x^{\star} \|$ and that
\[ \min \{ h_{x^k} (P_k), 0 \} \cdot \tfrac{\| \nabla f (x^k) \|^2 \| x^k -
   x^{\star} \|^2}{[f (x^k) - f (x^{\star})]^2} \tfrac{[f (x^k) - f
   (x^{\star})]^2}{\| x^k - x^{\star} \|^2} \leq \min \{ h_{x^k} (P_k), 0 \}
   \cdot \tfrac{[f (x^k) - f (x^{\star})]^2}{\| x^k - x^{\star} \|^2}  \]
since $\min \{ h_{x^k} (P_k), 0 \} \leq 0$.
Re-arranging the terms, we get
\begin{align}
  \tfrac{1}{f (x^{k + 1}) - f (x^{\star})} - \tfrac{1}{f (x^k) - f
  (x^{\star})} ={} & \tfrac{f (x^k) - f (x^{\star}) - [f (x^{k + 1}) - f
  (x^{\star})]}{[f (x^{k + 1}) - f (x^{\star})] [f (x^k) - f (x^{\star})]}
  \nonumber\\
  \geq{} & \tfrac{- \min \{ h_{x^k} (P_k), 0 \}  \tfrac{[f (x^k) - f
  (x^{\star})]^2}{\| x^k - x^{\star} \|^2}}{[f (x^{k + 1}) - f (x^{\star})] [f
  (x^k) - f (x^{\star})]} \label{eqn:pf-lem-6.1-4} \\
  ={} & \tfrac{- \min \{ h_{x^k} (P_k), 0 \}  [f (x^k) - f (x^{\star})]}{[f
  (x^{k + 1}) - f (x^{\star})] \| x^k - x^{\star} \|^2} \nonumber\\
  \geq{} & \tfrac{- \min \{ h_{x^k} (P_k), 0 \}}{\| x^k - x^{\star} \|^2} \geq -
  \tfrac{1}{\Delta^2} \min \{ h_{x^k} (P_k), 0 \}, \label{eqn:pf-lem-6.1-5}
\end{align}
where \eqref{eqn:pf-lem-6.1-4} plugs in \eqref{eqn:pf-lem-6.1-3}; \eqref{eqn:pf-lem-6.1-5} uses the fact that $f(x^k) \leq f(x^1)$ and that
\[\|x^k - x^\star\| \leq \max_{x \in \Lcal_{f(x^1)}} \|x - x^\star\| = \Delta.\]
Finally, we telescope the relation
\[  \tfrac{1}{f (x^{k + 1}) - f (x^{\star})} - \tfrac{1}{f (x^k) - f
  (x^{\star})} \geq -
  \tfrac{1}{\Delta^2} \min \{ h_{x^k} (P_k), 0 \}\]
from $k = 1$ to $K$ and get
\begin{align}
  \tfrac{1}{f (x^{K + 1}) - f (x^{\star})} - \tfrac{1}{f (x^1) - f
  (x^{\star})} ={} & \textstyle \sum_{k = 1}^K \tfrac{1}{f (x^{k + 1}) - f (x^{\star})} -
  \tfrac{1}{f (x^k) - f (x^{\star})} \nonumber\\
  \geq{} & - \tfrac{1}{\Delta^2} \textstyle \sum_{k = 1}^K \min \{ h_{x^k} (P_k), 0 \}
  \nonumber\\
  ={} & \tfrac{1}{\Delta^2} \textstyle \sum_{k = 1}^K \max \{ - h_{x^k} (P_k), 0 \} . \nonumber
\end{align}
Re-arranging the terms and using convexity of $\max\{\cdot, 0\}$,
\begin{align}
  f (x^{K + 1}) - f (x^{\star}) \leq & \tfrac{\Delta^2}{\sum_{k = 1}^K \max \{ -
  h_{x^k} (P_k), 0 \}}
  \leq  \tfrac{\Delta^2}{K} \tfrac{1}{\max \{ - \frac{1}{K} \sum_{k = 1}^K
  h_{x^k} (P_k), 0 \}} \nonumber
\end{align}

and this completes the proof.

\subsection{Proof of Proposition \ref{prop:hx-learnability}}
To show the Lipschitzness of $h_x$, it suffices to show the gradient is bounded.
Given $\nabla h_x (P) = \frac{\nabla f (x - P \nabla f (x)) \nabla f
(x)^{\top}}{\| \nabla f (x) \|^2},$ we deduce that
\begin{align}
  \| \nabla h_x (P) \|_F ={} & \tfrac{\| \nabla f (x - P \nabla f (x)) \nabla f
  (x)^{\top} \|_F}{\| \nabla f (x) \|^2} \nonumber\\
  ={} & \tfrac{\| \nabla f (x - P \nabla f (x)) \|}{\| \nabla f (x) \|} \label{eqn:prop-6.1-1}\\
  \leq{} & \tfrac{\| \nabla f (x - P \nabla f (x)) - \nabla f (x) \| + \|
  \nabla f (x) \|}{\| \nabla f (x) \|} \label{eqn:prop-6.1-2} \\
  \leq{} & \tfrac{L \| P \nabla f (x) \|}{\| \nabla f (x) \|} + 1 \label{eqn:prop-6.1-3} \\
  \leq{} & L \| P \| + 1 \leq L D + 1, \nonumber
\end{align}

where \eqref{eqn:prop-6.1-1} uses $\| a b^{\top} \|_F = \| a \| \cdummy \| b \|$ and \eqref{eqn:prop-6.1-3} applies
$L$-Lispschitzness of $\nabla f (x)$.

\subsection{Proof of Lemma \ref{lem:hx-learnability}}

The proof is again a direct application of the results in online convex
optimization. For any $P \in \mathcal{P}$, \eqref{eqn:pf-lem4.1-2} gives
\[ \textstyle \sum_{k = 1}^K h_{x^k} (P_k) - \sum_{k = 1}^K h_{x^k} (P) \leq \tfrac{1}{2
   \eta} \| P_1 - P \|_F^2 + \tfrac{\eta}{2} \sum_{k = 1}^K \| \nabla h_{x^k}
   (P_k) \|^2_F \leq \tfrac{2 D^2}{\eta} + \tfrac{\eta (LD + 1)^2}{2} K, \]
where the last inequality $\| P_1 - P \|_F \leq \| P_1 \|_F + \| P \|_F \leq 2 D$ and the
bounded gradient $\| \nabla h_{x^k} (P) \|_F \leq L (D + 1)$.
Taking $\eta$ to minimize the right-hand side completes the proof.

\subsection{Proof of Theorem \ref{thm:hx-trajectory-opt}}

By \Cref{lem:hx-learnability}, we have
\[ \textstyle \tfrac{1}{K} \sum_{k = 1}^K h_{x^k} (P_k) \leq \tfrac{1}{K} \sum_{k = 1}^K
   h_{x^k} (P) + \tfrac{\rho_K}{K} \]
for all $P \in \Pcal$
and plugging $- \tfrac{1}{K} \sum_{k = 1}^K h_{x^k} (P_k) \geq -
\theta_P^{\star} - \tfrac{\rho_K}{K}$ into \Cref{lem:hx-surrogate-measure} completes the
proof.

\subsection{Proof of Lemma \ref{lem:hx-hindsight}}

According to \ref{A3}, $L^{- 1} I \in \Pcal$ and descent lemma gives, for all $x \not \in \mathcal{X}^{\star}$, that
\[ h_x (L^{- 1} I) = \tfrac{f ( x - \frac{1}{L} \nabla f (x) ) -
   f (x)}{\| \nabla f (x) \|^2} \leq - \tfrac{1}{2 L} \]
and this completes the proof.

\subsection{Proof of Corollary \ref{coro:hx-globalconv}}

Using \Cref{lem:hx-hindsight} and \Cref{thm:hx-trajectory-opt}, $\theta_K^\star \leq - \gamma^\star$ and plugging the bound back into \Cref{thm:hx-trajectory-opt} completes the proof.

\end{document}